\documentclass[reqno]{amsart}

\usepackage{amsmath}
\usepackage{amssymb}
\usepackage[psamsfonts]{eucal}
\usepackage{fancyhdr}

\newtheorem{teo}{\sc Theorem}
\newtheorem{cor}{\sc Corollary}
\newtheorem{lemma}{\sc Lemma}
\newtheorem{defi}{\sc Definition}
\newtheorem{rem}{\sc Remark}

\newcommand{\adj}{\operatorname{Adj}}
\newtheorem{example}{Example}

\title[Relative asymptotics for orthogonal matrix polynomials]{Relative asymptotics for orthogonal matrix polynomials }

\author[A. Branquinho]{\sc A. Branquinho}
\address[A. Branquinho]{CMUC and Departamento de Matem\'atica, Universidade de Coimbra,
Apartado 3008, EC Universidade, 3001-454 Coimbra, Portugal.}
\email[A. Branquinho]{ajplb@mat.uc.pt}

\author[F. Marcell\'{a}n]{\sc F. Marcell\'{a}n}
\address[F. Marcell\'{a}n]{Departamento de Matem\'aticas,
Escuela Polit\'ecnica Superior, Universidad Carlos III de Madrid, Avenida de
la Universidad, 30, 28911 Legan\'es, Spain.}
\email[F. Marcell\'{a}n]{pacomarc@ing.uc3m.es}

\author[A. Mendes]{\sc A. Mendes}
\address[A. Mendes]{Departamento de Matem\'atica, Escola Superior de Tecnologia e Gest\~{a}o,
Instituto Polit\'ecnico de Leiria,  2411 - 901  Leiria - Portugal.}
\email[A. Mendes]{aimendes@estg.ipleiria.pt}

\thanks{{\hspace{-.45 cm}}{\it 2000 Mathemathics Subject Classification}. 33C45, 39B42. 
The work of the first (AB) and third (AM) authors has been supported by CMUC, Department of Ma\-the\-ma\-tics,
University of Coimbra. The work of the second author  has been supported by   Direcci\'on General de Investigaci\'on, Ministerio de Ciencia e Innovaci\'on of Spain, grant  MTM2009-12740-C03-01.}

\begin{document}

\begin{abstract}
In this paper we study sequences of matrix polynomials that sa\-tis\-fy a non-symmetric
recurrence relation. To study this kind of sequences we use a vector  interpretation of the matrix orthogonality.
In the context  of these sequences of matrix polynomials we introduce the concept of the generalized matrix Nevai class and we give the ratio asymptotics between two consecutive polynomials belonging to this class. We study the generalized matrix Chebyshev polynomials and  we  deduce  its  explicit expression as well as we show some illustrative examples.
The concept of a Dirac delta functional is introduced. We show how the vector model that includes a Dirac delta functional is a representation of a discrete Sobolev inner product. It also allows to reinterpret such  perturbations in the usual matrix Nevai class. Finally,  the relative asymptotics between a polynomial in the generalized matrix Nevai class and a polynomial that is orthogonal to a modification of the corresponding matrix measure by the addition of a Dirac delta functional is deduced.
\end{abstract}

\maketitle

\section{Introduction}  \label{sec:1}

In the last decade the asymptotic behavior of matrix orthonormal polynomials, the distribution of their zeros as well as  their connection with matrix quadrature formulas have paid a increasing attention by many researchers. A big effort was done in this direction in the framework of  the analytic theory of such polynomials by A. J. Dur\'{a}n, W. Van Assche and coworkers, among others (cf.~\cite{Dur93,DurWal95,DurDaneriVias,simon,sinapassche}).

In this work we study outer ratio asymptotics for matrix orthogonal polynomials belonging to a new class, the so  called {\it generalized matrix Nevai class,} and for these matrix orthogonal polynomials we obtain some new analytic results.
 \begin{defi}
Let $A$, $B$, and $C$ be matrices, with $A$  and $C$ non-singular matrices of dimension  $N\times N$. A {\it sequence of matrix polynomials} $\{V_m\}_{m \in {\mathbb{N}}}$ defined~by
 \begin{eqnarray}
 \label{ana3} zV_m(z) =A_{m} V_{m+1}(z) + B_m V_{m}(z) + C_m V_{m-1}(z), \quad m \geq 0,\end{eqnarray}
with $A_m$ a non-singular lower triangular matrix and $C_m$ a non-singular upper-triangular matrix, belongs to the {\it generalized matrix Nevai class}
$M(A,B,C)$~if  $$\displaystyle \lim_{m \to \infty} A_m =A, \quad \displaystyle \lim_{m \to \infty} B_m =B, \quad \mbox{and} \, \displaystyle \lim_{m \to \infty} C_m =C.$$
We say that a {\it matrix of measures} ${\mathcal{M}}$ belongs to the {\it generalized matrix Nevai class} $M(A,B,C)$ if some of the corresponding sequences of matrix orthogonal polynomials belongs to $M(A,B,C)$.
 \end{defi}

In order to to study the generalized matrix Nevai class we will recover a vector interpretation of the matrix orthogonality that was presented for the first time in~\cite{anamendes}. Let $(\mathbb{P}^{N})^*$ be the linear space of vector linear functionals defined on the linear space  $\mathbb{P}^{N}$ of vector polynomials  with complex coefficients.$(\mathbb{P}^{N})^*$ is said to be the {\it dual space.} A
\textit{vector of functionals} $
{\mathcal{U}}= [
 u^{1} \cdots \,
u^{N} ] ^{T}$ acting in $ \mathbb{P}^{N}$
over~$ \mathcal{M} _{N\times N}(\mathbb{C})$~is defined~by
 \begin{eqnarray*}
\mathcal{U}(\mathcal{P}):= (\mathcal{U}.\mathcal{P} ^{T} )
^{T}=\left[
 \begin{matrix}
\langle u^1, p_1 \rangle  & \cdots & \langle u^N, p_1 \rangle  \\
\vdots & \ddots & \vdots \\
\langle u^1, p_N \rangle  & \cdots & \langle u^N, p_N \rangle
 \end{matrix}
\right] \, ,
 \end{eqnarray*}
where \ ``$ . $'' \ means the symbolic product of the vectors
$ \mathcal{U}$ and $
\mathcal{P}^{T},$ where ${\mathcal{P}}^T= [p_1 \, \cdots \, p_N ]$, $p_i \in {\mathbb{P}}$, the standard linear space of polynomials with complex coefficients. The degree of ${\mathcal{P}}$ is given~by
 $
\operatorname{deg}({\mathcal{P}})= \lfloor (\underset{j=1, \ldots ,N}{\max } \{ \operatorname{deg}\text{
}p_{j} \} )/N \rfloor
,
 $
where $ \lfloor  \pmb{.}  \rfloor $
represents the integer part of a real number.

Given a polynomial $h$, with  $\deg h = N,$ the set
 \begin{eqnarray*}
\{1,x,\ldots,x^{N-1}, h(x), x h(x), \ldots, x^{N-1}h(x),h^2(x), xh^2(x), \ldots\}
\end{eqnarray*}
 is a basis for the linear space of polynomials,~${\mathbb{P}}$. Furthermore, $\{{{\mathcal{P}}_j}\}_{j \in {\mathbb{N}}},$
with ${\mathcal{P}}_j$ defined by $ {\mathcal{P}}_j(x)=(h(x))^j{\mathcal{P}}_0(x)$, where  ${\mathcal{P}}_0(x)=[1 \, x\, \cdots \, x^{N-1}]^T,$ is also a basis for the linear space of vector polynomials~${\mathbb{P}}^N.$

The vector of li\-ne\-ar functionals $(x^k {\mathcal{U}})$  acting in $ \mathbb{P}^{N}$
over~$ \mathcal{M} _{N\times N}(\mathbb{C})$ is defined~by
$$
 (x^k {\mathcal{U}})(\mathcal{P}):=((x^k {\mathcal{U}}).\mathcal{P}^T)^T={\mathcal{U}}(x^k{\mathcal{P}}) \, .
$$
So, with this definition and taking into account that $\{{{\mathcal{P}}_j}\}_{j \in {\mathbb{N}}}$ is a basis for the linear space of vector polynomials~${\mathbb{P}}^N$, the {\it j-$th$ moment} associated with the vector of linear functionals~$x^k {\mathcal{U}}$ is
given by $(x^k{\mathcal{U}})({\mathcal{P}}_j) = {\mathcal{U}}^k_j$.

The {\it Hankel matrices} associated with $\mathcal{U}$ are the matrices
 \begin{eqnarray*} 
 D_m =\left[%
 \begin{matrix}
  {\mathcal{U}}_0 & \cdots &  {\mathcal{U}}_m  \\
  \vdots & \ddots & \vdots \\
   {\mathcal{U}}_m & \cdots &  {\mathcal{U}}_{2m}
 \end{matrix}%
 \right],\,  m \in {\mathbb{N}},
 \end{eqnarray*}
where ${\mathcal{U}}_j$ is the  $j$-th moment associated with the vector of linear functionals ${\mathcal{U}}$.
${\mathcal{U}}$ is said to be {\it quasi-definite} if all the leading principal submatrices of $D_m, \, m \in {\mathbb{N}},$ are non-singular.

A {\it vector sequence of polynomials} $\{{\mathcal{B}}_m\}_{m \in {\mathbb{N}}}$, with degree of ${\mathcal{B}}_m$ equal to $m$, is said to be {\it left-orthogonal} with respect to the vector of linear functionals~${\mathcal{U}}$~if
$$(h^{k} {\mathcal{U}}) \left( {\mathcal{B}}_m \right)=
\Delta_m \delta_{k,m}, \, k=0,\ldots, m-1, \, m \in {\mathbb{N},}$$ where $\delta_{k,m}$ is the  Kronecker delta and $\Delta_m$ is a non-singular upper triangular matrix.

Notice that, we can always write $\mathcal{B}_{m}$ in the matrix form
 \begin{eqnarray*}
\mathcal{B}_{m}(x)=V_{m}(h(x))\mathcal{P}_{0}(x),
\end{eqnarray*}
where $ V_{m}$ is a$m$ degree~$N\times N$ matrix polynomial and $\mathcal{P}_{0}(x) = [
 1 \, x\, \cdots \, x^{N-1} ]^{T}.$

Similarly, a {\it sequence of matrix polynomials} $\{G_m\}_{m \in {\mathbb{N}}}$, with degree of $G_m$ equal to $m$, is said to be {\it right-orthogonal} with respect to the vector of linear functionals~${\mathcal{U}}$~if
 $$(G_m^T(h(x)) {\mathcal{U}}) \left( {\mathcal{P}}_j \right)=
\Theta_m\delta_{j,m},\, j=0,1, \ldots ,m-1, \,m \in \mathbb{N}$$ where $\Theta_m$ is a non-singular lower triangular matrix.

In~\cite{anamendes} necessary and sufficient conditions for the quasi-definiteness of~${\mathcal{U}}$, i.e., for the existence of a vector (matrix) sequence of polynomials left-orthogonal (right-orthogonal) with respect to the vector of linear functionals~$\mathcal{U}$ are obtained.

These sequences of polynomials satisfy non-symmetric three-term recurrence relations. So, if $\{{\mathcal{B}}_m\}_{m \in {\mathbb{N}}}$ is a vector sequence of polynomials left-orthogonal with respect to ${\mathcal{U}}$ and if $\{G_m\}_{m \in {\mathbb{N}}}$ is a sequence of matrix polynomials, then there exist sequences of numerical matrices $\{A_m\}_{m \in {\mathbb{N}}}$, $\{B_m\}_{m \in {\mathbb{N}}}$, and $\{C_m\}_{m \in {\mathbb{N}}}$, with $A_m$ a non-singular lower triangular matrix and $C_m$ a non-singular upper triangular matrix, such that
\begin{eqnarray} \label{rr1} h(x)
{\mathcal{B}}_m(x) = A_{m} {\mathcal{B}}_{m+1}(x) + B_m
{\mathcal{B}}_m(x) + C_m {\mathcal{B}}_{m-1}(x), \quad m\geq 1,
\end{eqnarray}
with
$ {\mathcal{B}}_{-1}(x) = 0_{1 \times N}$ and $ {\mathcal{B}}_{0}(x)={\mathcal{P}}_0(x) \, $,
where $ {\mathcal{P}}_0(x) = [1\, x\, \cdots \, x^{N-1} ]^T$ and
\begin{eqnarray} \label{rrrG}
z  G_n(z) =
G_{n-1}(z)A_{n-1} + G_n(z)B_n+ G_{n+1}(z)C_{n+1}, \ n\geq 1 \, ,
\end{eqnarray}
with $  G_{-1}(z)=0_{N \times N}$ and $G_0(z)={\mathcal{U}}({\mathcal{P}}_0)^{-1}$.

Notice  that these three-term recurrence relations completely characterize  each type of orthogonality.

Furthermore, right and left vector orthogonality are connected with right and left matrix orthogonality. Indeed, consider the {\it generalized Markov matrix function}~${\mathcal{F}}$ associated with~${\mathcal{U}}$  defined~by
 \begin{eqnarray*} 
{\mathcal{F}}(z) :=
{\mathcal{U}}_x \left( \frac{{\mathcal{P}}_{0}(x)}{z-h(x)}\right) =
\left[
 \begin{matrix}
\langle u^1_x , \frac{1}{z-h(x)}\rangle  & \cdots & \langle u^N_x , \frac{1}{z-h(x)}\rangle  \\
 \vdots & \ddots & \vdots \\
\langle u^1_x , \frac{x^{N-1}}{z-h(x)}\rangle & \cdots & \langle u^N_x , \frac{x^{N-1}}{z-h(x)}\rangle
\end{matrix}%
\right],
\end{eqnarray*}
with $z$ such that $|h(x)|<|z|$ for every  $x \in {\sf L}$ where
$  {\sf L} =  \cup_{j=1,\ldots,N} \,\mbox{supp} \, u^j_x \, $.
Here ${\mathcal{U}}_x$ represents the action of
${\mathcal{U}}$ on the variable $x$ and ${\mathcal{P}}_0(x) $ the same as before
In fact, see~\cite{anamendes}, the matrix sequence $ \{ G_{n}\}_{n \in {\mathbb{N}}}$ and the vector sequence $ \{ \mathcal{B}_{m} \}_{m \in {\mathbb{N}}}$ are bi-orthogonal with respect to $\mathcal{U}$, i.e.,
 \begin{eqnarray*}
 (   (  G_{n} (  h(x) )   ) ^{T}\mathcal{U}_{x} )
 (  \mathcal{B}_{m} )  =I_{N\times N} \, \delta_{n,m},\text{\hspace{
0.1in}}n,m\in \mathbb{N} \, ,
 \end{eqnarray*}
if and only if  matrix the sequences $ \{ G_{n} \}_{n \in {\mathbb{N}}}$ and $ \{ V_{m} \}_{m \in {\mathbb{N}}}$, where
\begin{eqnarray}\label{defiV}
\mathcal{B}_m (z) = V_m (h(z)) \mathcal{P}_0(z),
\end{eqnarray}
 are bi-orthogonal with respect to $\mathcal{F}$, i.e.
 \begin{eqnarray*}
\frac{1}{2\pi i}\int_{C}V_{m} (  z )  \mathcal{F} ( z )  G_{n}
( z )  dz= I_{N \times N} \, \delta _{n,m},\text{\hspace{0.1in}} n,m\in \mathbb{N},\
 \end{eqnarray*}
where $C$ is a closed path in $\{z \in {\mathbb{C}}: |z| > |h(x)|, x \in {\sf L} \}$.

The sequences of matrix polynomials $\{V_m\}_{m \in {\mathbb{N}}}$ and $\{G_m\}_{m \in {\mathbb{N}}}$ presented here are orthogonal with respect to a matrix of measures which is  not necessarily positive definite. These sequences satisfy the three-term recurrence relations~(\ref{ana3}) and~(\ref{rrrG}), respectively. On the other hand, these recurrence relations yield  a characterization of  right and left matrix orthogonality.


The sequences of matrix polynomials
$  \{
\mathcal{B}_{m}^{(1)}\}_{m \in {\mathbb{N}}}$ and $  \{G_{m}^{(1)}\}_{m \in {\mathbb{N}}}$ given~by
\begin{eqnarray*}
\mathcal{B}_{m}^{(1)}(z) & = & \mathcal{U}_{x}\left(
\frac{V_{m+1}(z)-V_{m+1}(h(x) )}{z-h(x)}\,\mathcal{
P}_{0}(x)\right), \\
 G^{(1)}_{m} (z) & = & \left[ \left( \frac{G_{m+1}^T
(z)- G_{m+1}^T (h(x))}{z-h(x)}\right){\mathcal{U}}_x
\right]({\mathcal{P}}_0(x)) ,
\end{eqnarray*}
are said to be the \textit{sequences of associated polynomials of the first kind for} $  \mathcal{U}$ and  $  \{
\mathcal{B}_{m} \}_{m \in {\mathbb{N}}}$ and $\{G_{m}\}_{m \in {\mathbb{N}}}$, respectively. Here $  \mathcal{U}_{x}$ represents the action of $  \mathcal{U}$ on the variable $  x$.

As a consequence of the definition of associated polynomials of the first kind for $\mathcal{B}_{m}$ and $G_m$, we get
\begin{gather}
\label{ana1} V_{m+1}(z) {\mathcal{F}}(z) - {\mathcal{B}}^{(1)}_m (z) = {\mathcal{U}}_x
\left( \frac{{\mathcal{B}}_{m+1}(x)}{z-h(x)} \right),\\
\label{ana2}{\mathcal{F}}(z) G_{m+1}(z) - G^{(1)}_{m}(z) = \left(G_{m+1}^T (h(x)){\mathcal{U}}_x\right)\left(
\frac{{\mathcal{P}}_0(x)}{z-h(x)}\right).
\end{gather}

The sequences of associated polynomials of the first kind
$\{{\mathcal{B}}^{(1)}_m\}_{m \in {\mathbb{N}}}$, $\{G^{(1)}_m\}_{m \in {\mathbb{N}}}$ satisfy respectively the three-term recurrence relations~\eqref{ana3} and~\eqref{rrrG} with initial conditions ${\mathcal{B}}^{(1)}_{-1}(z)=0_{N \times N}$, ${\mathcal{B}}^{(1)}_{0}(z)=A_0^{-1}$ and $G^{(1)}_{-1}(z)=0_{N \times N}$, $G^{(1)}_{0}(z)=C_1^{-1}$.

The sequences of matrix polynomials $\{V_m\}_{m \in {\mathbb{N}}}$ and $\{G_m\}_{m \in {\mathbb{N}}}$ satisfy a Chris\-toffel-Darboux type formula
\begin{eqnarray} \label{cdm1}
(x-z)\sum_{k=0}^m G_k(z)V_k(x)=G_m(z)A_m
V_{m+1}(x)-G_{m+1}(z) C_{m+1}V_m(x),
\end{eqnarray}
with $x, \, z \in {\mathbb{C}}$ and its confluent form
\begin{eqnarray} \label{cdm11}
\sum_{k=0}^m
G_k(x)V_k(x)=G_m(x)A_m V_{m+1}^{\prime}(x)-G_{m+1}(x)C_{m+1}
V_m^{\prime}(x),
\end{eqnarray} with $x \in {\mathbb{C}}$. The Christoffel-Darboux formula characterizes the matrix or\-tho\-go\-na\-lity~(cf. \cite{ab}) and allows  us to deduce the following result.
\begin{teo}[Liouville-Ostrogradski type formula]\label{liouville}Let $\{V_m\}_{m \in {\mathbb{N}}}$ and $\{G_m\}_{m \in {\mathbb{N}}}$ be the sequences of matrix  polynomials bi-orthogonal with respect to ${\mathcal{F}}$ with $\{V_m\}_{m \in {\mathbb{N}}}$ defined by~\eqref{defiV}. Let $\{{\mathcal{B}}_m^{(1)}\}_{m \in {\mathbb{N}}}$ and $\{G_m^{(1)}\}_{m \in {\mathbb{N}}}$ be, respectively, the sequences of associated matrix polynomials of the first kind  for  $\{{\mathcal{B}}_m\}_{m \in {\mathbb{N}}}$ and $\{G_m\}_{m \in {\mathbb{N}}}$.
Then,
 \begin{eqnarray}\label{liouv}
{\mathcal{B}}^{(1)}_{m}G_{m} - V_{m+1} G^{(1)}_{m-1} = A_{m}^{-1}
 \end{eqnarray}
where $A_m$ is the non-singular coefficient that appears in the recurrence relation~\eqref{rr1}.
\end{teo}
\begin{proof}
The sequences of polynomials $\{V_m\}_{m \in {\mathbb{N}}}$, $\{G_m\}_{m \in {\mathbb{N}}}$, $\{\mathcal{B}^{(1)}_{m-1}\}_{m \in {\mathbb{N}}}$, and $\{G_{m-1}^{(1)}\}_{m \in {\mathbb{N}}}$ satisfy, respectively, the recurrence relations~\eqref{ana3} and~\eqref{rrrG} with standard initial conditions.

To prove this result we proceed by induction.
For $m=0$ the result follows from the initial conditions.
We assume the formula
$${\mathcal{B}}^{(1)}_{p}G_{p} - V_{p+1} G^{(1)}_{p-1} = A_{p}^{-1},$$
is true for $p=1,\ldots,m-1$.
To prove that this relation is also valid for  $p=m$  we consider the following steps.
First, we use the recurrence relation in
 ${\mathcal{B}}^{(1)}_{m}$ and in $V_{m+1}$, i.e.,
 \begin{multline*}
 {\mathcal{B}}^{(1)}_{m}G_{m} - V_{m+1} G^{(1)}_{m-1} = A_m^{-1}(zI_{N \times N}-B_m) ({\mathcal{B}}^{(1)}_{m-1}G_m -V_m G_{m-1}^{(1)} )\\-
A_m^{-1}C_m ({\mathcal{B}}^{(1)}_{m-2}G_m - V_{m-1}G_{m-1}^{(1)} ).
 \end{multline*}
Second, we prove that ${\mathcal{B}}^{(1)}_{m-1}G_m -V_m G_{m-1}^{(1)}=0_{N \times N}.$

Multiplying by $G_m$ in the right side of~\eqref{ana1} and multiplying  by $V_m$ in the left side of~\eqref{ana2}, changing $m$ by $m-1$ in the relations~\eqref{ana1} and~\eqref{ana2}, and then subtracting these equations, we get
\begin{eqnarray*}
 {\mathcal{B}}^{(1)}_{m-1}(z)G_m(z) -V_m(z) G_{m-1}^{(1)}(z).
\end{eqnarray*}
Then,
\begin{multline*}
{\mathcal{B}}^{(1)}_{m-1}(z)G_m(z) -V_m(z) G_{m-1}^{(1)}(z)\\ =V_m(z) (G_m^T(h(x)){\mathcal{U}}_x)\left(\frac{{\mathcal{P}}_0(x)}{z-h(x)}\right) - {\mathcal{U}}_x\left(\frac{{\mathcal{B}}_m(x)}{z-h(x)}\right)G_m(z).
\end{multline*}
Adding and subtracting $(G_m^T(h(x)){\mathcal{U}}_x)\left(\frac{{\mathcal{B}}_m(x)}{z-h(x)}\right)$ in  the last relation and taking in consideration the left and right orthogonalities, the result follows.

Using the above  result  we have
 \begin{eqnarray} \label{sl}{\mathcal{B}}^{(1)}_{m}G_{m} - V_{m+1} G^{(1)}_{m-1}= -
A_m^{-1}C_m \left({\mathcal{B}}^{(1)}_{m-2}G_m - V_{m-1}G_{m-1}^{(1)}\right).
\end{eqnarray}

Again, using the recurrence relations for  $G_{m}$ and $G_{m-1}^{(1)}$
in ${\mathcal{B}}^{(1)}_{m-2}G_m - V_{m-1}G_{m-1}^{(1)}$ we get
 \begin{multline*}
 {\mathcal{B}}^{(1)}_{m-2}G_m - V_{m-1}G_{m-1}^{(1)} =
 ({\mathcal{B}}^{(1)}_{m-2}G_{m-1}- V_{m-1}G_{m-2}^{(1)})(zI_{N \times N}-B_{m-1})C_m^{-1} \\ +
(V_{m-1} G^{(1)}_{m-3}-{\mathcal{B}}^{(1)}_{m-2} G_{m-2}) A_{m-2}C_m^{-1}.
 \end{multline*}
Since  ${\mathcal{B}}^{(1)}_{m-2}G_{m-1} -V_{m-1} G_{m-2}^{(1)}=0_{N \times N}$, we get
 \begin{eqnarray*}
{\mathcal{B}}^{(1)}_{m-2}G_m - V_{m-1}G_{m-1}^{(1)} = (V_{m-1} G^{(1)}_{m-3}
-{\mathcal{B}}^{(1)}_{m-2} G_{m-2}) A_{m-2}C_m^{-1}.
 \end{eqnarray*}
Using this relation in~\eqref{sl} we obtain
\begin{eqnarray*}{\mathcal{B}}^{(1)}_{m}G_{m} - V_{m+1} G^{(1)}_{m-1}=-
A_m^{-1}C_m (V_{m-1} G^{(1)}_{m-3}-{\mathcal{B}}^{(1)}_{m-2} G_{m-2}) A_{m-2}C_m^{-1}.
\end{eqnarray*}

According to the hypothesis of induction the result follows.
\end{proof}

Let ${\mathcal{U}}$ be a quasi-definite vector of linear functionals and let $\{V_m\}_{m \in {\mathbb{N}}}$ and $\{G_m\}_{m \in {\mathbb{N}}}$ be
sequences of bi-orthogonal polynomials. We denote the {\it kernel polynomial}~by
$$K_m(x,y)=\sum_{k=0}^{m-1} G_k(h(y))V_k(h(x)).$$ Notice that, even though we don't have $K_m(x,y)= K^T_m (y,x)$ like in the matrix symmetric case,  the reproducing pro\-per\-ty for the kernel holds.
\begin{teo}
Let ${\mathcal{U}}$ be a quasi-definite vector of linear functionals,
$\{{\mathcal{B}}_m\}_{m \in {\mathbb{N}}}$ and $\{G_m\}_{m \in {\mathbb{N}}}$ be,  respectively, the left  vector and right matrix orthogonal polynomials with respect to ${\mathcal{U}}$.  Given a vector polynomial $\pi \, \in {\mathbb{P}}^N$ of degree $m$, i.e.,
 \begin{eqnarray}\label{xx}
 \pi(x)=\sum_{k=0}^m \beta_k^m {\mathcal{B}}_k(x),
\quad \beta_k^m \in {\mathcal{M}}_{N \times N}({\mathbb{C}}),
\end{eqnarray}
then $\pi(x)= (K_{m+1}^T(x,z)
{\mathcal{U}}_z)(\pi(z)){\mathcal{P}}_0(x).$
\end{teo}
\begin{proof}
 From~\eqref{xx}, using the bi-orthogonality and taking into account
${\mathcal{L}}_n=G_n^T{\mathcal{U}}$, we get
$\beta_k^m = {\mathcal{L}}_k (\pi(z))=(G_k^T(h(z)){\mathcal{U}}_z)(\pi(z)).
$ Then, $$\pi(x)=\sum_{k=0}^m (G_k^T(h(z)){\mathcal{U}}_z)(\pi(z)){\mathcal{B}}_k(x).$$
But,
 $$
(G_k^T(h(z)){\mathcal{U}}_z)(\pi(z)){\mathcal{B}}_k(x)=((G_k(h(z)) V_k(h(x)))^T {\mathcal{U}}_z))(\pi(z)){\mathcal{P}}_0(x).
 $$
 Hence,
 the result follows.
\end{proof}
Let $J$ be the block matrix defined~by
 \begin{eqnarray} \label{J}
 J=
  \left[%
 \begin{matrix}
  B_0 & A_0 & 0 &  \\
  C_1 & B_1 & A_1 & \ddots \\
  0 & C_2 & B_2 & \ddots \\
    & & \ddots& \ddots
\end{matrix}%
\right]\end{eqnarray}  that is known in the literature as {\it $N$-Jacobi matrix}. When the matrix polynomials satisfy a
symmetric recurrence relation, it was proved in~\cite{DurLR} that the zeros of the $m$-th
orthogonal polynomial are the eigenvalues of $J_m$. This result can be generalized for sequences of orthogonal polynomials  satisfying
non-symmetric recurrence relations. Thus, for $m \in {\mathbb{N}}$,
the zeros of the matrix polynomials $G_{m}$ and $V_{m}$ are the zeros of the polynomial $\det (t I_{mN\times mN}- J_{m})$, with the same
multiplicity, where $I_{mN \times mN}$ is the identity matrix of
dimension $mN\times mN$ and $J_m$ is the truncated $N$-block Jacobi matrix
of dimension $mN\times mN$.

Also, for a matrix of measures $W$ and for any polynomial $V$ with degree less than  or equal to $2m-1$
the following quadrature formula holds
 \begin{eqnarray*}
 \int
V(h(x))dW(h(x)) = \sum_{k=1}^s V(x_{m,k})\Gamma_{m,k}
 \end{eqnarray*}
where  $x_{m,k},$ $k=1,\ldots,s,$
are the zeros of the matrix polynomial $V_m$, in general, complex numbers as well as  $s\leq mN$,  and
$\Gamma_{m,k}$ are the matrices
 \begin{eqnarray*}
\Gamma_{m,k}= \frac{l_k}{(\det\, (V_m(x)))^{(l_k)}(x_{m,k})}
(\adj\,(V_m(x)))^{(l_k -1)}(x_{m,k})
{\mathcal{B}}^{(1)}_{m-1}(x_{m,k}),
 \end{eqnarray*}
 for $k=1,\ldots,s,$
where $l_k$ is the multiplicity of the zero $x_{m,k}$.

Using the above quadrature formula, in~\cite{anamendes} we have obtained the following asym\-pto\-tic result:
$$\lim_{m \to \infty} V^{-1}_m(z) {\mathcal{B}}^{(1)}_{m-1} (z) = {\mathcal{F}}(z)$$
locally uniformly in ${\mathbb{C}}\setminus \Gamma$, where
$ \Gamma=\cap_{N\geq 0} M_N$,  $ M_N=\overline{\cup_{n\geq N}\{\mbox{zeros of }\,
V_n \}} \, .$

The structure of the manuscript is as follows.
 In section~\ref{sec:2}, sequences of matrix orthogonal polynomials belonging to the ge\-ne\-ra\-li\-zed matrix Nevai class are studied. Furthermore,  the outer ratio asymptotics  of two consecutive polynomials belonging to this class is obtained. We also study the generalized matrix Chebyshev polynomials and we present their explicit formulas as well as the  corresponding  generalized Markov function. The example presented in this section does not belong to the cases studied by A. J. Dur\'{a}n in~\cite{Dur99} and also can't be converted in these cases by the method presented by H. Dette and coworkers (see~\cite{Dettealunos}).

In section~\ref{sec:3}, we introduce a modification of a vector of linear functionals  by  adding a Dirac delta. This yields  a reinterpretation and an extension of a perturbation in the usual matrix Nevai class (see~\cite{nevai,MarcellanYakhlefPinar1,MarcellanYakhlefPinar2}). On the other hand, the meaning of  any discrete Sobolev inner product in vector terms is clarified and this is clearly related with sequences of polynomials satisfying higher order recurrence relations (cf.~\cite{anamendes,Dur93,DurWal95,MarcellanZagorodnyuk,Zagorodnyuk2,Zagorodnyuk3,Zagorodnyuk4,Zagorodnyuk5}).
We also find  necessary and sufficient conditions for the quasi-definiteness of the modified functional.
The quasi-definiteness conditions obtained by the authors in~\cite{AlfaroMarcellanRezolaRonveaux,MarcellanAlfaroRezola,MarcellanPerezPinar2} for some special examples coincide with our general result.
To conclude this section, we describe the generalized Markov function associated to a modification by a Dirac delta functional.

Finally, in section~\ref{sec:4}  we present the relative asymptotics between a sequence of matrix  polynomials that belongs to the {\it generalized matrix Nevai class} and a sequence of matrix polynomials  that is orthogonal with respect to a modification by a Dirac delta functional of this class. This result generalizes  those obtained by F. Marcell\'{a}n and coworkers (see~\cite{MarcellanYakhlefPinar1}).

\section{Generalized matrix Nevai class}  \label{sec:2}

Let ${\mathcal{M}}$ be a matrix of measures in the generalized matrix Nevai class $M(A,B,C)$. Notice that  ${\mathcal{M}}$  can belong to several Nevai classes because  of the non-uniqueness of the corresponding sequences of orthogonal polynomials.

If $A$ and $C$ are non-singular matrices we can introduce the sequence of matrix polynomials $\{U_m^{A,B,C}\}_{m \in {\mathbb{N}}}$ defined by the recurrence formula
 \begin{eqnarray}\label{ana4}
 z U_m^{A,B,C}(z) = A U_{m+1}^{A,B,C}(z) + BU_m^{A,B,C}(z) +C U_{m-1}^{A,B,C}(z), \quad m\geq 1,
 \end{eqnarray}
with initial conditions $U_0^{A,B,C}(z)=I_{N \times N}$ and $U_{-1}^{A,B,C}(z)=0_{N \times N}$. According to an extension of the matrix Favard's theorem (see~\cite{anamendes}) this sequence is orthogonal with respect to a matrix of measures ${\mathcal{M}}_{A,B,C}$ that is  not necessarily positive definite. This sequence of matrix polynomials  is said to be the {\it sequence of matrix generalized second kind Chebyshev polynomials}.

The continued fraction associated with~\eqref{ana4} or, equivalently, the generalized Markov function ${\mathcal{F}}_{A,B,C}$ is given~by
 $$
\displaystyle
{\mathcal{F}}_{A,B,C} (z)
 =
\frac{1}{(zI_{N \times N}-B) -A\frac{1}{\displaystyle z I_{N \times N}-B-A\frac{1}{zI_{N\times N}-B-\cdots}C}C}
 $$
where $1/X$ denotes the inverse of the matrix $X$.

Matrix continued fractions of different types were studied by many authors (cf.~\cite{AptekarevNikishin,Fair,sorokin,Zygmunt}). In~\cite{Zygmunt} the reader can find  a detailed study about matrix continued fractions and matrix Chebyshev polynomials, where the author emphasyzes  how continued fractions are used to develop the notion of matrix Chebyshev polynomials in some symmetric cases.

Now, let us consider a sequence of vector  polynomials $\{\mathcal{B}_m\}_{m \in {\mathbb{N}}}$ left-orthogonal with res\-pect to the vector of linear functionals ${\mathcal{U}}$ satisfying the recurrence relation
 \begin{eqnarray*}
 h(z) \mathcal{B}_m(z) = A \mathcal{B}_{m+1}(z) + B\mathcal{B}_m(z) +C \mathcal{B}_{m-1}(z), \quad m\geq 1,
 \end{eqnarray*}
with initial conditions
${\mathcal{B}}_{-1}(z)=0_{N \times 1} $ and ${\mathcal{B}}_{0}(z)={\mathcal{P}}_0(z)$,
where $A$ and $C$ are non-singular matrices. It is straightforward to prove that the sequence of matrix polynomials $\{V_m\}_{m \in {\mathbb{N}}},$  defined by
$\mathcal{B}_m(z)=V_m (h(z)){\mathcal{P}}_0(z),$
and the sequence $\{{\mathcal{B}}^{(1)}_m\}_{m \in {\mathbb{N}}}$, of associated polynomials of the first kind for  ${\mathcal{U}}$ and $\{\mathcal{B}_m\}_{m \in {\mathbb{N}}}$, satisfies  the same recurrence relation with the following initial conditions $V_{-1}(z)=0_{N \times N}$, $V_0(z)=I_{N \times N}$, ${\mathcal{B}}^{(1)}_{-1}(z)=0_{N \times N},$ and ${\mathcal{B}}^{(1)}_0(z)=A^{-1}$, respectively. These conditions  are said to be the standard ones.

Rewriting the recurrence equations for $\{V_m\}_{m \in {\mathbb{N}}}$ and $\{{\mathcal{B}}^{(1)}_{m}\}_{m \in {\mathbb{N}}}$ in a blocks matrix form, we have
\begin{eqnarray*} \label{eq2}
 \left[
       \begin{matrix}
         A V_{m+1} & A {\mathcal{B}}^{(1)}_{m} \\
         V_m &  {\mathcal{B}}^{(1)}_{m-1}
       \end{matrix}
     \right] = \left[
                 \begin{matrix}
                   zI_{N \times N} - B & -C \\
                   I_{N \times N} & 0_{N \times N}
                 \end{matrix}
               \right]\left[
                        \begin{matrix}
                          V_m & {\mathcal{B}}^{(1)}_{m-1} \\
                          V_{m-1} & {\mathcal{B}}^{(1)}_{m-2}
                        \end{matrix}
                      \right].
\end{eqnarray*}
Since the matrix $A$ is non-singular  the last equation is equivalent to
\begin{eqnarray*} 
 \left[\begin{matrix}
                          V_{m+1} & {\mathcal{B}}^{(1)}_{m} \\
                          V_{m} & {\mathcal{B}}^{(1)}_{m-1}
                        \end{matrix}
                      \right]=\left[
                 \begin{matrix}
                   A^{-1}(zI_{N \times N} - B) & -A^{-1}C \\
                   I_{N \times N} & 0_{N \times N}
                 \end{matrix}
               \right]\left[
                        \begin{matrix}
                          V_m & {\mathcal{B}}^{(1)}_{m-1} \\
                          V_{m-1} & {\mathcal{B}}^{(1)}_{m-2}
                        \end{matrix}
                      \right].
 \end{eqnarray*}
Writing the last equation as
$L_{m}=TL_{m-1}$
where
$$L_{m}=\left[\begin{matrix}
                          V_{m+1} & {\mathcal{B}}^{(1)}_{m} \\
                          V_{m} & {\mathcal{B}}^{(1)}_{m-1}
                        \end{matrix}
                      \right]\quad \mbox{and} \quad T=\left[\begin{matrix}
  A^{-1}(zI_{N \times N} - B) & -A^{-1}C \\
                   I_{N \times N} & 0_{N \times N}
                 \end{matrix}
               \right],$$
we have
$L_{m}=T^{m}L_0$
with  $$L_0=\left[\begin{matrix}
                   A^{-1}(zI_{N \times N} - B) & A^{-1} \\
                   I_{N \times N} & 0_{N \times N}
                 \end{matrix}
               \right].$$
 For some particular choices of $A$, $B,$ and $C$ the matrix $T^m$ has  the spectral decomposition
$T^m=SD^mS^{-1}$
where  $D$  is a diagonal.

Using this decomposition we can determine $L_{m}$ and then, obtain $V_m$ and ${\mathcal{B}}^{(1)}_{m-1}$.
By a straightforward calculation we obtain $V_m^{-1}\,{\mathcal{B}}^{(1)}_{m-1}$ and then taking the limit when $m \to \infty$, we get
the generalized Markov function ${\mathcal{F}}_{A,B,C}$.

As a sake of example, in the case of matrix polynomials of dimension $2 \times 2$, if we consider the same problem but with different initial conditions, for instance,
$\widehat{V}_{0}(z)=P$, $\widehat{V}_1(z)=M+Qz$, $\widehat{{\mathcal{B}}}^{(1)}_{-1}(z)=0_{2 \times 2},$ and $\widehat{{\mathcal{B}}}^{(1)}_0(z)=Q$,~with
$$P=\left[%
\begin{matrix}
  p_{11} & p_{12} \\
 p_{21} & p_{22}
\end{matrix}%
\right], \,\, M=\left[%
\begin{matrix}
m_{11} & m_{12} \\
  m_{21} & m_{22}
\end{matrix}%
\right] \quad \mbox{and} \quad Q=\left[
                               \begin{matrix}
                                q_{11} & q_{12} \\
                                 q_{21} & q_{22}
                               \end{matrix}
                             \right]
,$$ we can relate the generalized Markov function associated with these new initial conditions, $\widehat{{\mathcal{F}}}_{A,B,C},$ with ${\mathcal{F}}_{A,B,C}$, in the following way: \\
\phantom{ol} -- The block matrix $\widehat{L}_m$ is given~by
\begin{eqnarray*}\widehat{L}_m=\widehat{T}^m \widehat{L}_0=T^m \widehat{L}_0=T^m L_0 L_0^{-1}\widehat{L}_0=L_m L_0^{-1}\widehat{L}_0.
\end{eqnarray*}
Using  this relation we get,
\begin{eqnarray*}
\widehat{V}_m(z)&=&V_m(z)P+{\mathcal{B}}^{(1)}_{m-1}(z)\left[A M +BP + (A Q-P)z\right]\\
\widehat{{\mathcal{B}}}^{(1)}_{m-1}(z)&=&{\mathcal{B}}^{(1)}_{m-1}(z)A Q\\
\widehat{{\mathcal{F}}}_{A,B,C}(z)&=&P^{-1} {\mathcal{F}}_{A,B,C}(z) A Q + M^{-1} Q + P^{-1} B^{-1} A Q + (I-P^{-1} A Q)z.
\end{eqnarray*}
In the next example, we illustrate how can we determine the generalized Markov function associated
with generalized second kind Chebyshev polynomials for a particular choice of the recurrence coefficients with standard initial conditions.
\begin{example}
Let us consider   $A$, $B,$ and $C$ as $$A= \left[%
\begin{matrix}
  1 & 0 \\
 0 & 1
\end{matrix}%
\right], \,\, B=\left[%
\begin{matrix}
-1 & 0 \\
  1 & -1
\end{matrix}%
\right], \quad \mbox{and} \quad C=\left[
                               \begin{matrix}
                                -1 & 0 \\
                                 0 & 1
                               \end{matrix}
                             \right]
.$$
The matrices $T$ and $L_0$ are
$$T=\left[\begin{matrix}
                                 1+z & 0 & -1 & 0\\
                                 1 & z+1 & 0 & 1 \\
                                 1 & 0 & 0 & 0 \\
                                 0 & 1 & 0 & 0
                               \end{matrix}
                             \right]\quad \mbox{and} \quad L_0=\left[
    \begin{matrix}
       1+z & 0 & 1 & 0 \\
       -1 & 1+z & 0 & 1 \\
       1 & 0 & 0 & 0 \\
       0 & 1 & 0 & 0
    \end{matrix}
  \right].$$
Matrix $T$ have the following eigenvalues
\begin{eqnarray*}
\lambda_1=\frac{1}{2}(1+z-\sqrt{(z+1)^2-4}),\, \lambda_2=\frac{1}{2}(1+z+\sqrt{(z+1)^2-4}),\\
\lambda_3=\frac{1}{2}(1+z-\sqrt{(z+1)^2+4}),\, \lambda_4=\frac{1}{2}(1+z+\sqrt{(z+1)^2+4})\end{eqnarray*}
and the corresponding eigenvectors are
$$v_1=\left[-2,\frac{1}{2}(1+z-\sqrt{(z+1)^2+4}),
\frac{4}{-1-z+\sqrt{(z+1)^2+4}},1\right],$$
$$v_2=\left[-2,\frac{1}{2}(1+z+\sqrt{(z+1)^2+4}),
-\frac{4}{1+z+\sqrt{(z+1)^2+4}},1\right],$$
$$v_3=\left[0,\frac{1}{2}(1+z-\sqrt{(z+1)^2-4}),0,1\right],$$
$$v_4=\left[0,\frac{1}{2}(1+z+\sqrt{(z+1)^2-4}),0,1\right].$$
Then, the matrix $T^m$ has the following spectral decomposition
$$T^m=SD^mS^{-1},$$
where $D =\operatorname{diagonal}\,[ \lambda_1,  \lambda_2,  \lambda_3, \lambda_4],$
for $m \in {\mathbb{N}},$ and  $S=[v_1 |\, v_2 |\, v_3 |\, v_4]$, where $\lambda_i$ and $v_i$, for $i=1,\ldots,4$ are, respectively, the eigenvalues and the eigenvectors of the matrix $D$.
Thus  we can determine $L_{m+1}$ and then, obtain $V_m$ and ${\mathcal{B}}^{(1)}_{m-1}$ as follows
\begin{multline*}
V_m(z)=-\frac{1}{2}E_{m+2}(z)\left[
              \begin{matrix}
                0 & 0 \\
                1+z & 0
              \end{matrix}
            \right]+E_{m+1}(z)\left[
              \begin{matrix}
                1+z & 0 \\
               -\frac{1}{2}& 0
              \end{matrix}
            \right]+E_{m}(z)\left[
              \begin{matrix}
                1 & 0 \\
                0 & 0
              \end{matrix}
            \right] \\ + F_{m+1}(z) \left[
              \begin{matrix}
                0 & 0 \\
               \frac{1}{2}(2+z)z & 1+z
              \end{matrix}
            \right] - F_{m}(z) \left[
              \begin{matrix}
                0 & 0 \\
               1+z& 1
              \end{matrix}
            \right],
            \end{multline*}
\begin{gather*}
{\mathcal{B}}^{(1)}_{m-1}(z)
 =
E_{m+2}(z)\left[
              \begin{matrix}
                0 & 0 \\
                -\frac{1}{2} & 0
              \end{matrix}
            \right]
            +
E_{m+1}(z)\left[
              \begin{matrix}
                1 & 0 \\
                0 & 0
              \end{matrix}
            \right]
            +
F_{m+1}(z)\left[
              \begin{matrix}
                0 & 0 \\
               \frac{1}{2}(1+z) & 1
              \end{matrix}
            \right] ,
\end{gather*}
where $$
E_m(z)=2^{-m}\frac{(1+z+\sqrt{(z+1)^2 +4})^m-(1+z-\sqrt{(z+1)^2 +4})^m}{\sqrt{(1+z)^2+4}},$$
$$F_m(z)=2^{-m}\frac{(1+z+\sqrt{(z+1)^2 -4})^m-(1+z-\sqrt{(z+1)^2 -4})^m}{\sqrt{(1+z)^2-4}}.$$
By a straightforward calculation we obtain $V_m^{-1}\,{\mathcal{B}}^{(1)}_{m-1}.$ Taking the limit when{\linebreak}$m \to \infty$, we get
 $$
{\mathcal{F}}_{A,B,C}(z)=\left[
                          \begin{matrix}
                            \frac{2}{1+z+\sqrt{(1+z)^2+4}}  & 0 \\
                            \frac{4+(1+z-\sqrt{(1+z)^2+4})(1+z-\sqrt{(1+z)^2-4}-\sqrt{(1+z)^2+4})}{(1+z+\sqrt{(1+z)^2-4})(1+z+\sqrt{(1+z)^2+4})}  &\frac{2}{1+z+\sqrt{(1+z)^2-4}}
                          \end{matrix}
                        \right] .
 $$

\end{example}
This example inspired us to prove the ratio asymptotics  between two con\-se\-cu\-ti\-ve polynomials in the generalized matrix Nevai class. Indeed,  without loss of generality, we consider  a sequence of matrix polynomials $\{V_m\}_{m \in {\mathbb{N}}}$ defined by~\eqref{ana3} such that $V_0(z)=I_{N \times N}$. But first, we need to introduce the following result.
 \begin{lemma}\label{kill}
Let $\{V_m\}_{m \in {\mathbb{N}}}$ be a sequence of matrix orthogonal polynomials in the generalized matrix Nevai class $M(A,B,C)$.
Then, there exists a positive constant~$M$, which does not depend on $m$, such that their zeros $x_{m,k}$ are  contained in a disk  $D=\{z \in {\mathbb{C}}:|z| < M\}$.
 \end{lemma}
 \begin{proof}
Let us consider the $N$-block Jacobi matrix, $J$ (see~\eqref{J}), associated with the recurrence relation~\eqref{ana3}.
Remember also, that the zeros of $V_{m}$ are the eigenvalues of~$J_m$
where $J_m$ is the truncated matrix of $J$, with dimension $mN \times mN$.

Taking into account that the sequences $\{A_m\}_{m \in {\mathbb{N}}}$, $\{B_m\}_{m \in {\mathbb{N}}}$, $\{C_m\}_{m \in {\mathbb{N}}}$ converge, and using the Gershgorin disk theorem for the location of eigenvalues, it follows that there exists $M>0$ such that if $x_{m,k}$ is a zero of $V_m$ then $x_{m,k} \in D$ where $D=\{z \in {\mathbb{C}}:|z| < M\}$.  So, $\Gamma$ defined~by
$$\Gamma=\cap_{N\geq 0} M_N,\quad M_N=\overline{\cup_{m\geq N}\{\mbox{zeros of}\,\,
V_m \}},$$
is contained in $D$ and $\mbox{supp}(W)\subset \Gamma \subset D$.
 \end{proof}
 \begin{teo}\label{maintheorem}
Let $\{V_m\}_{m \in {\mathbb{N}}}$ be a sequence of matrix polynomials left-ortho\-gonal with respect to the matrix of measures ${\mathcal{M}}$ and satisfying the three-term recurrence relation~\eqref{ana3}. Assume that $\displaystyle \lim_{m \to \infty} A_m =A$, $\displaystyle\lim_{m \to \infty} B_m =B,$ and $\displaystyle \lim_{m \to \infty} C_m =C$ with $A$ and $C$ non-singular matrices. Then,
 \begin{eqnarray}\label{assrelativa}
\lim_{m \to \infty} V_{m-1}(z) V_m^{-1}(z) A_{m-1}^{-1}= {\mathcal{F}}_{A,B,C}(z), \quad z \in {\mathbb{C}}\setminus \Gamma,\end{eqnarray}
where ${\mathcal{F}}_{A,B,C}$ is the Markov transform of the  matrix of measures for the generalized second kind Chebyshev polynomials. Moreover, the convergence is locally  uniformly  for compact subsets of ${\mathbb{C}}\setminus \Gamma$, where $\Gamma=\cap_{N\geq 0} M_N,\,
M_N=\overline{\cup_{m\geq N}Z_m},$ and  $Z_m$ is the set of the zeros of $V_m$.
 \end{teo}
 \begin{proof}
First, we consider the sequence of discrete matrix measures{\linebreak}$\{\mu_m\}_{m \in {\mathbb{N}}}$ defined~by
$$\mu_m = \sum_{k=1}^s V_{m-1}(h(y_{m,k})) \Gamma_{m,k} G_{m-1}(h(y_{m,k})) \,\delta_{y_{m,k}} , \quad m \geq 0,$$
where $y_{m,k}$ are complex numbers such that $h(y_{m,k})=x_{m,k},$  with $x_{m,k}$, $k=1, \ldots, s,$ the zeros of the polynomial $V_m$, and the matrix  $\Gamma_{m,k}$ is given~by
 \begin{eqnarray*}
 \Gamma_{m,k}= \frac{l_k(\adj\,(V_m(x)))^{(l_k
-1)}(x_{m,k}) {\mathcal{B}}^{(1)}_{m-1}(x_{m,k})}{(\det\,
(V_m(x)))^{(l_k)}(x_{m,k})}, \, k=1,\ldots,s,\end{eqnarray*}
 $l_k$ being the multiplicity of the zero $x_{m,k}$, $l_k \leq N,$  and $\{{\mathcal{B}}^{(1)}_{m-1}\}_{m \in {\mathbb{N}}}$ the sequence of associated polynomials of the first kind for $\{{\mathcal{B}}_m\}_{m \in {\mathbb{N}}}$ and ${\mathcal{U}}$.
 Notice that  $\Gamma_{m,k}$ is the weight in the quadrature formula. Hence, it follows that
$$\int d\mu_m (h(x))=I_{N \times N}, \quad \mbox{for}\quad m \geq 0.$$
The decomposition, (cf.~\cite{gant}),
$$V_{m-1}(z)V_m^{-1}(z) = \sum_{k=0}^s C_{m,k}\frac{1}{z-x_{m,k}},$$
with $$C_{m,k}=\frac{l_kV_{m-1}(x_{m,k})(\adj\,(V_m(t)))^{(l_k
-1)}(x_{m,k})}{(\det\,
(V_m(t)))^{(l_k)}(x_{m,k})}$$
is always possible even though the zeros of $V_m$ are complex or have multiplicity greater than one (see~\cite{Cberg,anamendes,Dur96}).
Then, we have
$$C_{m,k}A_{m-1}^{-1}= \frac{l_kV_{m-1}(x_{m,k})(\adj\,(V_m(t)))^{(l_k
-1)}(x_{m,k})A_{m-1}^{-1}}{(\det\,
(V_m(t)))^{(l_k)}(x_{m,k})}.$$
Applying the generalized Liouville-Ostrogradski formula~\eqref{liouv} and taking into account that for every $\textbf{b},$ a zero of a matrix polynomial $V_m,$  (see~\cite{Dur95})
 \begin{eqnarray*} 
V_m(\textbf{b})\left(\adj\, (V_m(t))\right)^{(p-1)}(\textbf{b})=\left(\adj\, (V_m(t))\right)^{(p-1)}(\textbf{b})V_m(\textbf{b})=0_{N \times N},
 \end{eqnarray*}
we obtain
 \begin{eqnarray*}
 C_{m,k}A_{m-1}^{-1}= V_{m-1}(h(y_{m,k})) \Gamma_{m,k}G_{m-1}(h(y_{m,k})).
 \end{eqnarray*}
 From the definition of the matrix of measures $\mu_m$ we get
$$V_{m-1}(z)V_m^{-1}(z)A_{m-1}^{-1}=\int \frac{d\mu_m(h(x))}{z-h(x)}, \quad z \in {\mathbb{C}} \setminus \Gamma.$$

Let us consider the generalized Chebyshev matrix polynomials of second kind $\{U_m^{A,B,C}\}_{m \in {\mathbb{N}}}$ defined by~\eqref{ana4}. We can prove by induction that
 \begin{eqnarray} \label{kill1}
 \lim_{m \to \infty} \int U_l^{A,B,C} (h(t)) d\mu_m(h(t)) = \left\{
                                                                 \begin{array}{ll}
                                                                   I_{N \times N}, & \hbox{for $l=0$,} \\
                                                                   0_{N \times N}, & \hbox{for $l\neq 0$.}
                                                                 \end{array}
                                                               \right.
 \end{eqnarray}
To do it, we just  use the same technicalities as in~\cite{Dur99}.

We are now ready to prove
$$
 \lim_{m \to \infty} \int \frac{d\mu_m(h(x))}{z-h(x)}= {\mathcal{F}}_{A,B,C}(z), \quad z \in {\mathbb{C}} \setminus \Gamma.
$$
If not, we can find a complex number $z \in {\mathbb{C}} \setminus \Gamma$, an increasing sequence of nonnegative integers $\{n_{l}\}_{l \in {\mathbb{N}}}$, and a positive constant $C$ such that
 \begin{eqnarray} \label{kill2}
\left\|\int \frac{d\mu_{n_l}(h(x))}{z-h(x)} - {\mathcal{F}}_{A,B,C}(z)\right\|_2 \geq C >0, \,\, l\geq0,
 \end{eqnarray}
where  $\|\, . \,\|_2$ denotes the {\it spectral norm of a matrix}, i.e.,
$$\|A\|_2=\mbox{max}\{\sqrt \lambda: \, \lambda \, \mbox{is a eigenvalue of}\, A^* A\}.$$
Since $\{\mu_m\}_{m \in {\mathbb{N}}}$ is a sequence of matrices of measures with support contained in a disk $D$ (see lemma~\ref{kill}) and taking into account that $  \int d\mu_m = I_{N \times N}$, by using the Banach-Alaoglu theorem, we can obtain  a subsequence $\{r_l\}_{l \in {\mathbb{N}}}$ from $\{n_l\}_{l \in {\mathbb{N}}}$, defined on a curve $\gamma_M$ contained in the disk $D$, with the same $k$-th moments of the vector of linear functionals ${\mathcal{U}}$, for $k \leq 2 r_l -1$, such that
$$\lim_{l \to \infty} \int_{\gamma_M} f(h(x))d\mu_{r_l}(h(x))= \frac{1}{2\pi i} \int_{\gamma_M} f(h(z)){\mathcal{U}}_x\left(\frac{{\mathcal{P}}_0(x)}{z-h(x)}\right)dz,$$
for any continuous matrix function  $f$ defined in $D$.

Hence,  taking $f(h(x))=U^{A,B,C}(h(x))$, we have
 \begin{eqnarray*}
 \lim_{l\to \infty} \int_{\gamma_M}  U^{A,B,C}_l(h(x)) d\mu_{r_l}(h(x))= \frac{1}{2\pi i} \int_{\gamma_M} U^{A,B,C}_l(h(z)){\mathcal{U}}_x\left(\frac{{\mathcal{P}}_0(x)}{z-h(x)}\right).
 \end{eqnarray*}
 From~\eqref{kill1} we have
\begin{eqnarray*}
 \frac{1}{2\pi i} \int_{\gamma_M} U^{A,B,C}_l(h(z)){\mathcal{U}}_x\left(\frac{{\mathcal{P}}_0(x)}{z-h(x)}\right)= \left\{
                                                                 \begin{array}{ll}
                                                                   I_{N \times N}, & \hbox{for $l=0$,} \\
                                                                   0_{N \times N}, & \hbox{for $l\neq 0$.}
                                                                 \end{array}
                                                               \right.
 \end{eqnarray*}
But, the sequence of matrix polynomials $\{U_m^{A,B,C}\}_{m \in {\mathbb{N}}}$ is orthogonal with respect to ${\mathcal{F}}_{A,B,C}$. Since $\{U_m^{A,B,C}\}_{m \in {\mathbb{N}}}$ is a basis of the linear space of matrix polynomials we get that~\eqref{kill2} is not possible.
Each entry of the matrix
$ \int \frac{d\mu_{m}(h(x))}{z-h(x)}$ is uniformly bounded on compact sets of
${\mathbb{C}}\setminus \Gamma$. Then, according to the Stieltjes-Vitali theorem, we get the uniform  convergence.
\end{proof}

Notice that we have analogous results of Lemma~\ref{kill} and Theorem~\ref{maintheorem} for the sequence of matrix polynomials $\{G_m\}_{m \in {\mathbb{N}}}$.
 \begin{cor}
Under the hypothesis of Theorem \ref{maintheorem} we have that
\begin{equation}\label{convderivadas}
\lim_{m \to \infty}(V_{m-1}(z) V_m^{-1}(z))^{(k)}A_{m-1}^{-1}={\mathcal{F}}_{A,B,C}^{(k)}(z) \end{equation}
on compact subsets of ${\mathbb{C}}\backslash \Gamma$, for $k=1,2, \ldots\, .$
 \end{cor}
The locally  uniformly  convergence in~\eqref{convderivadas} means that every entry of the left hand-side of~\eqref{convderivadas} is locally uniformly convergent to its corresponding entry in the right hand-side of~\eqref{convderivadas}.

\section{Delta functionals}  \label{sec:3}

 In this section we deal with  a vector of linear functionals that results of a modification by a Dirac delta functional and we will illustrate it with a nice application related to Sobolev inner products. Thus, we can reinterpret these inner products in a vectorial form and it is  a  motivation for the study of these modifications.

In a generic way, we can say that a {\it Dirac delta functional} is a vector of linear functionals where the functional components are linear combinations of Dirac deltas and their derivatives in a finite set of points in the real line.

If the polynomial $h$ of fixed degree $N$ consider in the previous sections is such~that
 \begin{eqnarray} \label{h} h(x)=\prod_{j=1}^M
(x-c_j)^{M_j+1},
 \end{eqnarray}
where $M_j+1$ is the multiplicity of each $c_j
\, \in \,{\mathbb{N}}$ as a zero of $h,$ then we can define a new vector of linear functionals as the result of  a modification by the addition of  a Dirac  delta functional with respect to $h$ as follows.
 \begin{defi} The vector of linear functionals $\widetilde{\mathcal{U}}$ defined~by
 \begin{eqnarray} \label{utilde} \widetilde{{\mathcal{U}}}={\mathcal{U}}+ \Lambda \, \pmb{\delta},
 \end{eqnarray}
 with $\Lambda$ a numerical matrix of dimension $N\times N$ where
 \begin{eqnarray*}
\pmb{\delta} = [\delta_{c_1}\,
\delta^{\prime}_{c_1}\, \cdots \, \delta^{(M_1)}_{c_1}\, \delta_{c_2}
\, \delta^{\prime}_{c_2}\, \cdots \,\delta^{(M_2)}_{c_2} \, \cdots \,
\delta_{c_M}\, \delta^{\prime}_{c_M} \,\cdots \, \delta^{(M_M)}_{c_M}]^T,
 \end{eqnarray*}
where  $N=M + \sum_{j=1}^M M_j$, is called {\it  vector of linear functionals modified by a Dirac delta functional associated with  $h$}.
\end{defi}


First, as a motivation for the study of these mo\-di\-fi\-cations, we will consider the following example. It gives a vectorial reinterpretation of a Sobolev inner product. It is important to refer that with the same technics presented in this example, a vector  reinterpretation for any general discrete Sobolev inner product holds.
 \begin{example} [see~\cite{MarcellanPerezPinar2}]
Let us consider the discrete Sobolev inner product
 \begin{equation}\label{sobolev1}
 \langle f,g
\rangle_S:=\int_{I} f g d\mu + \lambda f^{\prime}(0)g^{\prime}(0), \quad
\mbox{where} \quad \lambda \in {\mathbb{R}}^+ \, .
 \end{equation}
To establish the parallelism between vector orthogonality and Sobolev inner products, let $\{\widetilde{p}_n\}_{n \in {\mathbb{N}}}$ be a sequence of scalar polynomials orthonormal with respect to the inner product~\eqref{sobolev1}, i.e.
 \begin{eqnarray}\label{cosob1}
\langle \widetilde{p}_n, x^k \rangle_S &=& 0,\,
k=0,\ldots,n-1,\\\nonumber \langle \widetilde{p}_n, x^n \rangle_S
&\neq& 0, \, n \in {\mathbb{N}}.
 \end{eqnarray}
Notice that the multiplication by the polynomial $h(x)=x^2$ yields a symmetric operator with respect to the above  Sobolev inner product, i.e.,
$$\langle x^2f,g \rangle_S=\langle f,x^2g\rangle_S, \, \forall f,\, g \in {\mathbb{P}},$$
and, as a consequence, the corresponding  sequence of orthonormal polynomials
$\{\widetilde{p}_n\}_{n \in {\mathbb{N}}}$ satisfies, for $n \geq 0$, a five-term recurrence relation
 \begin{eqnarray}
 \label{sobo5termos}
x^2 \widetilde{p}_{n+1} = c_{n+3,2}
\widetilde{p}_{n+3} + c_{n+2,1}
\widetilde{p}_{n+2}+c_{n+1,0}\widetilde{p}_{n+1}+ c_{n+1,1}
\widetilde{p}_{n}+c_{n+1,2} \widetilde{p}_{n-1} \, .
 \end{eqnarray}
On the other hand, taking into account the vectorial approach  given  in~\cite{anamendes}, the recurrence relation~\eqref{sobo5termos} yields  the following vector expression
 \begin{eqnarray*}
  x^2
\widetilde{{\mathcal{B}}}_m(x) = \widetilde{A}_{m+1}
\widetilde{{\mathcal{B}}}_{m+1}(x) + \widetilde{B}_m
{\mathcal{B}}_m(x) + \widetilde{A}^{T}_m {\mathcal{B}}_{m-1}(x), \,
m \geq 0,
 \end{eqnarray*}
 where
 $\widetilde{{\mathcal{B}}}_m (x)= [
\widetilde{p}_{2m}(x) \, \widetilde{p}_{2m+1}(x) ]^T$,
 \begin{eqnarray*} \widetilde{A}_m =  \left[
 \begin{matrix}
  c_{2m,2} & 0 \\
  c_{2m,1} & c_{2m+1,2}
 \end{matrix}%
 \right], \quad \mbox{and} \quad  \widetilde{B}_m = \left[
 \begin{matrix}
  c_{2m,0} & c_{2m+1,1}\\
  c_{2m+1,1} & c_{2m+1,0}
 \end{matrix}%
 \right].
 \end{eqnarray*}
Notice that $\widetilde{B}_m =\widetilde{B}^T_m$ as well as  we are dealing with a symmetric case.
Then, by Favard's type theorem (see~\cite{anamendes}), there exists a vector of linear functionals
 $\widetilde{{\mathcal{U}}} = [\widetilde{u}^1 \,\, \, \widetilde{u}^2 ]^T$
such that  the  sequence of vector polynomials
$\{\widetilde{{\mathcal{B}}}_m\}_{m \in {\mathbb{N}}}$ is left-ortho\-gonal with respect to $\widetilde{{\mathcal{U}}}$ and
  \begin{eqnarray} \label{cosob2}(x^{2k}
\widetilde{{\mathcal{U}}}) \left( \widetilde{{\mathcal{B}}
}_m\right)&=& 0_{2 \times 2},\,\, k=0,1, \ldots ,m-1,\\
(x^{2m} \widetilde{{\mathcal{U}}}) \left(
\widetilde{{\mathcal{B}}}_m \right)&=& \widetilde{\Delta}_m, \,\, m
\in \mathbb{N},\nonumber \end{eqnarray} where $\widetilde{\Delta}_m$ is a non-singular upper triangular matrix.

 From  the definition of a vector of linear functionals, the previous orthogonality conditions are given explicitly, for all $k=0,\ldots,m-1 $ and $m \in {\mathbb{N}}$~by
 $$
 (x^{2k} \widetilde{{\mathcal{U}}}) \left( \widetilde{{\mathcal{B}}}_m \right)= \left[
 \begin{matrix}
  \langle \widetilde{u}^1, x^{2k} \widetilde{p}_{2m}\rangle & \langle \widetilde{u}^2, x^{2k} \widetilde{p}_{2m}\rangle\\
  \langle \widetilde{u}^1, x^{2k} \widetilde{p}_{2m+1}\rangle & \langle \widetilde{u}^2, x^{2k} \widetilde{p}_{2m+1}\rangle
 \end{matrix} \right] = \left[
 \begin{matrix}
  0 & 0\\
  0 & 0
 \end{matrix}
 \right] ,
 $$
 $$
 (x^{2m} \widetilde{{\mathcal{U}}}) \left( \widetilde{{\mathcal{B}}}_m \right)= \left[
 \begin{matrix}
  \langle \widetilde{u}^1, x^{2m} \widetilde{p}_{2m}\rangle & \langle \widetilde{u}^2, x^{2m} \widetilde{p}_{2m}\rangle\\
  \langle \widetilde{u}^1, x^{2m}\widetilde{ p}_{2m+1}\rangle & \langle \widetilde{u}^2, x^{2m} \widetilde{p}_{2m+1}\rangle
 \end{matrix}
 \right] = \left[
 \begin{matrix}
  \bullet & \bullet\\
  0 & \bullet
 \end{matrix}
 \right] \, . $$

 The vector of linear functionals $\widetilde{{\mathcal{U}}}$ represents a Sobolev inner product like~\eqref{sobolev1} only if it is a modification by a Dirac delta functional, as we will describe in the sequel. In fact, the vector of linear functionals $\widetilde{{\mathcal{U}}}$ has the representation
 \begin{eqnarray} \label{lena}
 \widetilde{{\mathcal{U}}}={\mathcal{U}}+
\Lambda \, \pmb{\delta}
 \end{eqnarray} where
 $${\mathcal{U}}=\left[u \,\,
xu\right]^T, \quad
\Lambda= { \left[%
\begin{matrix}
   0 & 0   \\
   0 & - \lambda
\end{matrix}%
\right]}, \quad \pmb{\delta} =\left[\delta_0 \,\, \delta_0^{\prime}
\right]^T,
 $$
and where $u$ is a linear functional on the linear space of scalar polynomials ${\mathbb{P}}$ such~that
$u(p(x))=\int_I p(x)\, d\mu(x).$ Here $\mu$ is the weight function that appears in the Sobolev discrete inner product~\eqref{sobolev1}.

 To illustrate that the orthogonality conditions~\eqref{cosob1} and~\eqref{cosob2} are equivalent, first, take $k=0$ in~\eqref{cosob2},
\begin{eqnarray*}
 \widetilde{{\mathcal{U}}} \left( \widetilde{{\mathcal{B}}}_m
 \right)=\left[
 \begin{matrix}
 \langle 1,\widetilde{p}_{2m} \rangle_S & \langle x,\widetilde{p}_{2m} \rangle_S\\
 \langle 1,\widetilde{p}_{2m+1} \rangle_S& \langle x,\widetilde{p}_{2m+1} \rangle_S
 \end{matrix}
 \right]=
 \left[ \begin{matrix}
  0 & 0\\
  0 & 0
 \end{matrix}
 \right].
\end{eqnarray*}
 Furthermore, for all $k=1, \ldots, m, \, m \in {\mathbb{N}}$
\begin{eqnarray*}
(x^{2k}\widetilde{{\mathcal{U}}})\left( \widetilde{{\mathcal{B}}}_m
\right)&=&\left[ \begin{matrix}
 \langle x^{2k},\widetilde{p}_{2m} \rangle_S & \langle x^{2k+1},\widetilde{p}_{2m} \rangle_S\\
 \langle x^{2k},\widetilde{p}_{2m+1} \rangle_S& \langle x^{2k+1},\widetilde{p}_{2m+1} \rangle_S
\end{matrix}\right] =
\left[ \begin{matrix}
  \neq 0 & \bullet \\
  0 & \neq 0
\end{matrix}
 \right]\delta_{k,m} \, .
 \end{eqnarray*}
Thus, these  relations show us that the orthogonality conditions~\eqref{cosob1} and~\eqref{cosob2} are equivalent and it means that the discrete Sobolev inner product~\eqref{sobolev1} can be represented in the vectorial form by~\eqref{lena}.
\end{example}

Motivated by this example the first question naturally imposed consists to know when $\widetilde{{\mathcal{U}}}$ defined by~\eqref{utilde} is quasi-definite, i.e., when there exists a  sequence of vector  polynomials $\{\widetilde{{\mathcal{B}}}_m \}_{m \in {\mathbb{N}}}$ left-ortho\-gonal with respect to $\widetilde{{\mathcal{U}}}$. Then, the next step is to obtain necessary and sufficient conditions for the quasi-definiteness of $\widetilde{{\mathcal{U}}}.$

In the sequel, we denote by $\{{\mathcal{B}}_m\}_{m \in {\mathbb{N}}}$ the vector sequence of  polynomials left-ortho\-gonal with respect to  ${\mathcal{U}}$  and $\{\widetilde{{\mathcal{B}}}_m\}_{m \in {\mathbb{N}}}$ the sequence of polynomials asso\-cia\-ted with  $\widetilde{{\mathcal{U}}}$ defined by~\eqref{utilde}, i.e., the sequence of polynomials satisfying
 \begin{eqnarray}
\label{associada}\widetilde{{\mathcal{U}}} \left(
\widetilde{{\mathcal{B}}}_m \right)= 0_{N \times N},\, m \geq 1,
\quad \mbox{and} \quad  \widetilde{{\mathcal{U}}} \left(
\widetilde{{\mathcal{B}}}_0 \right) \,\, \mbox{is a non-singular matrix}.
\end{eqnarray}

As above, we will denote by $\{\widetilde{V}_m\}_{m \in {\mathbb{N}}}$ the sequence of matrix polynomials defined~by $$\widetilde{{\mathcal{B}}}_{m}(x)=\widetilde{V}_m(h(x)){\mathcal{P}}_0(x),$$
with ${\mathcal{P}}_0(x) = [1\, x\, \cdots\, x^{N-1} ]^T$.

Before the statement of conditions that give us quasi-definiteness for $\widetilde{{\mathcal{U}}}$, we need the following auxiliary result.
\begin{lemma}
Let $\{{\mathcal{B}}_m\}_{m \in {\mathbb{N}}}$ be a vector sequence of polynomials left-ortho\-gonal with res\-pect to the vector of linear functionals  ${\mathcal{U}}$ and $\{\widetilde{{\mathcal{B}}}_m\}_{m \in {\mathbb{N}}}$ be a vector sequence of polynomials associated with  $\widetilde{{\mathcal{U}}}$ defined by~\eqref{utilde} and verifying~\eqref{associada}. Then, the following statements  hold
\begin{eqnarray*}
 \nonumber
 (h^k\widetilde{{\mathcal{U}}})({\mathcal{B}}_m) &=& (h^k{\mathcal{U}})({\mathcal{B}}_m)=0_{N \times N}, \, k \geq 1,\\
 (h^k\widetilde{{\mathcal{U}}})(\widetilde{{\mathcal{B}}}_m)
&=& (h^k{\mathcal{U}})(\widetilde{{\mathcal{B}}}_m),\, k \geq 1, \, m \in {\mathbb{N}}.
 \end{eqnarray*}
\end{lemma}
\begin{proof}
The proof of this result is straightforward taking into account that  for any vector polynomial ${\mathcal{P}}\, \in \, {\mathbb{P}}^N$,  we
have
$$(h^k\pmb{\delta}) ({\mathcal{P}})=0_{N \times N}, \quad \mbox{for}\quad k\geq 1,$$
 with
 \begin{eqnarray*}
 \pmb{\delta} = [\delta_{c_1}\,
\delta^{\prime}_{c_1}\, \cdots \, \delta^{(M_1)}_{c_1}\, \delta_{c_2}
\, \delta^{\prime}_{c_2}\, \cdots \,\delta^{(M_2)}_{c_2} \, \cdots \,
\delta_{c_M}\, \delta^{\prime}_{c_M} \,\cdots \, \delta^{(M_M)}_{c_M}]^T,
\end{eqnarray*}
and the result follows.
\end{proof}
%
%

 In the literature about discrete Sobolev orthogonal polynomials (see for example~\cite{AlfaroMarcellanRezolaRonveaux,MarcellanAlfaroRezola,BavinckMeijer,EvansETal,MarcellanPerezPinar,MarcellanPerezPinar2,MarcellanRonveaux}) the reader could find the same conditions that we will achieve in the next result. Another interesting work related with this topic is~\cite{MarcellanYakhlefPinar1}. There, the authors studied the outer relative asymptotics  between a matrix polynomial belonging to a matrix Nevai class and a matrix polynomial that is orthogonal to a perturbation of such a matrix of measures in the matrix Nevai class. There the existence of matrix orthogonal polynomials with respect to perturbations in the matrix Nevai class is assumed.
The next theorem gives  necessary and sufficient conditions for the existence of such matrix polynomials in a more general case.
\begin{teo}
 \label{regularidade}Let ${\mathcal{U}}$ be a quasi-definite vector of linear functionals, $\{{\mathcal{B}}_m\}_{m \in {\mathbb{N}}}$ be the vector sequence of polynomials left-ortho\-gonal with respect to ${\mathcal{U}}$ and
$\{\widetilde{{\mathcal{B}}}_m\}_{m \in {\mathbb{N}}}$ be the associated vector sequence of polynomials, verifying~\eqref{associada}, to the vector linear functional $\widetilde{{\mathcal{U}}}$ defined by~\eqref{utilde}.
The vector of linear functionals $\widetilde{{\mathcal{U}}}$ is quasi-definite if and only~if
 \begin{eqnarray*}
I_{N \times N} + \pmb{\delta}_z ({\mathcal{P}}_0(z)) \Lambda^T
K_{m+1}(c_l,c_l), \, m \geq 0,
 \end{eqnarray*}
is a non-singular matrix for all
$c_l$ zero of $h$ defined by~\eqref{h}.
 \end{teo}
 \begin{proof}
For  the vector sequences $\{{\mathcal{B}}_{m}\}_{m \in {\mathbb{N}}}$ and
$\{\widetilde{{\mathcal{B}}}_{m}\}_{m \in {\mathbb{N}}}$, we can consider the sequences of matrix polynomials $\{V_m\}_{m \in {\mathbb{N}}}$ and $\{\widetilde{V}_m\}_{m \in {\mathbb{N}}}$ such that
$${\mathcal{B}}_{m}(z)=V_m(h(z)){\mathcal{P}}_0(z) \quad \mbox{and} \quad \widetilde{{\mathcal{B}}}_{m}(z)=\widetilde{V}_m(h(z)){\mathcal{P}}_0(z),$$
with ${\mathcal{P}}_0(z) = [1\,z\, \cdots\, z^{N-1} ]^T$.
Using the reproducing property for the kernel we get  $$V_m(h(x))= (K^T_{m+1}(x,z) {\mathcal{U}}_z)({\mathcal{B}}_m(z)).$$
Similarly, we can write
\begin{eqnarray*} \nonumber
\widetilde{V}_m(h(x)) =
(K_{m+1}^T(x,z) \widetilde{{\mathcal{U}}}_z )
(\widetilde{{\mathcal{B}}}_m (z))- \sum_{j=0}^m (G_j^T(h(z))\Lambda
\, \pmb{\delta}_z) (\widetilde{{\mathcal{B}}}_m (z))V_j(h(x)),
\end{eqnarray*}
where $\pmb{\delta}_z $ is the Dirac  delta functional acting on $z$.
 Notice that in the right hand side of the above identity
 \begin{eqnarray*}
(G^T_j(h(z))\Lambda \,\pmb{\delta}_z )({\widetilde{V}}_m
(h(z)){\mathcal{P}}_0(z)) =(G^T_j(h(c_l))\Lambda \,\pmb{\delta}_z )(\widetilde{V}_m (h(z)){\mathcal{P}}_0(z)),
\end{eqnarray*}
holds for any zero  $c_l$ of the polynomial $h$, independently of its multiplicity. From now on, we write $h(c_l)=0$.

Again, since
\begin{eqnarray*}
\pmb{\delta}_z (\widetilde{V}_m
(h(z)) {\mathcal{P}}_0(z))=\widetilde{V}_m(0)\pmb{\delta}_z ({\mathcal{P}}_0(z)),
\end{eqnarray*}
 then
 \begin{eqnarray*}\sum_{j=0}^m
(G_j^T(h(z))\Lambda \, \pmb{\delta}_z)(\widetilde{{\mathcal{B}}}_m
(z))V_j(h(x))=  \widetilde{V}_m(0)\sum_{j=0}^m
\pmb{\delta}_z ({\mathcal{P}}_0(z)) \Lambda^T G_j(0)V_j(h(x)).
\end{eqnarray*}

Now, to analyze
$(K_{m+1}^T(x,z)\widetilde{{\mathcal{U}}}_z ) (\widetilde{{\mathcal{B}}}_m (z))$ we
use the definition of the kernel $K_m(x,z)$ $=\sum_{j=0}^{m-1} G_j(h(z))V_j(h(x))$ and consider the following representation for $G_m$,
$G_m (h(x))=\sum_{k=0}^m \beta_k^m (h(x))^k$,  $\beta_k^m \in
{\mathcal{M}}_{N \times N}({\mathbb{C}}).$ Then, we have
 \begin{eqnarray*}
 (K_{m+1}^T(x,z) \widetilde{{\mathcal{U}}}_z )
(\widetilde{{\mathcal{B}}}_m (z))=\sum_{j=0}^m [(\sum_{k=0}^j
(h(x))^k\widetilde{{\mathcal{U}}}_z)( \widetilde{{\mathcal{B}}}_m
(z)) \beta_k^j ]V_j(h(x)).
\end{eqnarray*}
Hence,
\begin{eqnarray*}
(K_{m+1}^T(x,z) \widetilde{{\mathcal{U}}}_z )
(\widetilde{{\mathcal{B}}}_m (z)) = \widetilde{\Delta}_m \beta_m^m
V_m(h(x))
\end{eqnarray*}
if and only if
$\{\widetilde{{\mathcal{B}}}_m \}_{m \in {\mathbb{N}}}$
is left-ortho\-gonal with respect to
$\widetilde{{{\mathcal{U}}}}_z$, satisfying the or\-tho\-go\-na\-li\-ty conditions
\begin{eqnarray*}
(h^{k} \widetilde{{\mathcal{U}}}) \left(
\widetilde{{\mathcal{B}}}_m \right)&=& \widetilde{\Delta}_m\delta_{k,m}, \,
k=0,1, \ldots ,m,\, m \in \mathbb{N},
\end{eqnarray*}
where $\widetilde{\Delta}_m$ is a non-singular matrix.
Thus,
\begin{eqnarray} \label{fr}
\widetilde{V}_m(h(x))=D_m
V_m(h(x))-\widetilde{V}_m(0)\pmb{\delta}_z({\mathcal{P}}_0(z))\Lambda^T
K_{m+1}(x,c_l)
 \end{eqnarray}
where $D_m=\widetilde{\Delta}_m \beta_m^m$ is a non-singular matrix. Taking $x=c_l$ we get
$$\widetilde{V}_m(0)(I_{N \times N} + \pmb{\delta}_z ({\mathcal{P}}_0(z))\Lambda^T K_{m+1}(c_l,c_l) ) =  D_m V_m(0).$$
$\{\widetilde{V}_m\}_{m \in {\mathbb{N}}}$, is completely determined by the data if and only~if
$$I_{N \times N} + \pmb{\delta}_z ({\mathcal{P}}_0(z))\Lambda^T
K_{m+1}(c_l,c_l), m\geq 0,$$ is a non-singular matrix for any zero  $c_l$  of the polynomial $h$, independently of its multiplicity.
\end{proof}

If we apply this result to the previous example, then it is straightforward  to deduce  that these  conditions  are the  same as those obtained  in~\cite{MarcellanPerezPinar2} for the existence of a sequence of orthogonal polynomials with respect to the Sobolev inner product~\eqref{sobolev1}.

In the sequel  $h$ is a polynomial of fixed degree $N$ defined by~\eqref{h} and sequences of matrix polynomials $\{V_m\}_{m \in {\mathbb{N}}}$ and  $\{G_m\}_{m \in {\mathbb{N}}}$ are given by $V_m(z)=\sum_{j=0}^m \alpha_j^m z^j$ and $G_m(z)=\sum_{j=0}^m \beta_j^m z^j$. These sequences are bi-orthogonal with respect  to the generalized Markov function ${\mathcal{F}}$ and satisfy the recurrence relations~\eqref{ana3} and~\eqref{rrrG}, respectively.
\begin{teo}
 Let ${\mathcal{U}}$ be a quasi-definite vector functional, $\widetilde{{\mathcal{U}}}$ quasi-definite vector functional defined by~\eqref{utilde}, and $h$ be a  polynomial of fixed degree $N$. Let
${\mathcal{F}}$ and $\widetilde{{\mathcal{F}}}$ be the generalized Markov functions associated with ${\mathcal{U}}$ and
$\widetilde{{\mathcal{U}}}$, respectively. Then,
\begin{eqnarray}
\label{ftilde}z\widetilde{{\mathcal{F}}}(z)=z{\mathcal{F}}(z)+
\pmb{\delta}_x ({\mathcal{P}}_0(x))\Lambda^T,
\end{eqnarray}
with
${\mathcal{P}}_0(x) = [1\, x\, \cdots \, x^{N-1} ]^T$ and
$$\pmb{\delta} = [\delta_{c_1}\,
\delta^{\prime}_{c_1}\, \cdots \, \delta^{(M_1)}_{c_1}\, \delta_{c_2}
\, \delta^{\prime}_{c_2}\, \cdots \,\delta^{(M_2)}_{c_2} \, \cdots \,
\delta_{c_M}\, \delta^{\prime}_{c_M} \,\cdots \, \delta^{(M_M)}_{c_M}]^T,$$
with  $N=M + \sum_{j=1}^M M_j$.
 \end{teo}
 \begin{proof}
 From the definition of $\widetilde{{\mathcal{F}}}$ we have
 \begin{eqnarray*}
z\widetilde{{\mathcal{F}}}(z) &=& z\sum_{k=0}^\infty
\frac{\widetilde{{\mathcal{U}}}_x(h(x))^k{\mathcal{P}}_{0}(x)}{z^{k+1}}\\&=& z\sum_{k=0}^\infty
 \frac{{\mathcal{U}}_x ((h(x))^k{\mathcal{P}}_{0}(x))+ (\Lambda
\pmb{\delta}_x )((h(x))^k{\mathcal{P}}_{0}(x))}{z^{k+1}} \\&=&z {\mathcal{F}}(z)+
(\Lambda \pmb{\delta}_x ) \left( {\mathcal{P}}_{0}(x)\right)
\end{eqnarray*}
as we wanted to show.
\end{proof}

On the other hand, we can also consider the sequences $\{\widetilde{V}_m\}_{m \in {\mathbb{N}}}$ and $\{\widetilde{G}_m\}_{m \in {\mathbb{N}}}$ where $\widetilde{V}_m(z)=\sum_{j=0}^m \widetilde{\alpha}_j^m z^j$ and $\widetilde{G}_m(z)=\sum_{j=0}^m \widetilde{\beta}_j^m z^j$ are bi-ortho\-gonal sequences with respect to the functional $\widetilde{{\mathcal{F}}}$ defined by~\eqref{ftilde}.


The next result states the relation between  $\{V_m\}_{m \in {\mathbb{N}}}$ and  $\{\widetilde{V}_m\}_{m \in {\mathbb{N}}}.$
 \begin{teo} Let $\{V_m\}_{m \in {\mathbb{N}}}$ be a sequence of matrix polynomials left-ortho\-gonal with respect to the generalized Markov function ${\mathcal{F}}$ and satisfying the three-term recurrence relation
$$h(x)V_m(z)=A_m V_{m+1}(z) + B_m V_m(z)+ C_{m}V_{m-1}(z), \quad m\geq 1.$$
Let $\{\widetilde{V}_m\}_{m \in {\mathbb{N}}}$ be the sequence of matrix polynomials left-ortho\-gonal with respect to $\widetilde{{\mathcal{F}}}$ defined by~\eqref{ftilde}. Then, the sequence of matrix polynomials{\linebreak}
$\{\widetilde{V}_m\}_{m \in {\mathbb{N}}}$ is defined~by
 \begin{eqnarray}
 \label{lixo}
h(x)\widetilde{V}_m(z) = \alpha^{1}_{m+1} V_{m+1}(z)
+\alpha^{2}_{m} V_m(z) + \alpha^{3}_{m-1} V_{m-1}(z) \, .
 \end{eqnarray}
with $\alpha^{1}_{m+1}$, $\alpha^{2}_{m}$,
and $\alpha^{3}_{m-1}$ such that
 \begin{eqnarray*} \alpha^{1}_{m+1}&=&
D_m A_m -L_m \beta_0^m A_m
\\ \alpha^{2}_{m}&=& D_m B_m +L_m
\beta_0^{m+1} C_{m+1} \\ \alpha^{3}_{m-1}&=& D_m C_m ,
 \end{eqnarray*}
 where  $L_m=D_m V_m(0)[I_{N
\times N}+ \pmb{\delta}_z ({\mathcal{P}}_0(z))\Lambda^T
K_{m+1}(c_l,c_l)]^{-1}\pmb{\delta}_z ({\mathcal{P}}_0(z))\Lambda^T$ and $D_m =(\widetilde{\beta}_m^m)^{-1} \beta_m^m$.
\end{teo}
\begin{proof}
 From the proof of Theorem~\ref{regularidade} we have
\begin{multline*}\widetilde{V}_m(h(x))= D_m V_m(h(x))\\- D_m V_m(0)[I_{N \times N}+ \pmb{\delta}_z ({\mathcal{P}}_0(z))\Lambda^T
K_{m+1}(c_l,c_l)]^{-1}\pmb{\delta}_z ({\mathcal{P}}_0(z))\Lambda^T
K_{m+1}(x,c_l). \end{multline*}
For the sake of simplicity, we denote by $L_m$ the matrix
$$L_m=D_m V_m(0)[I_{N \times N}+ \pmb{\delta}_z ({\mathcal{P}}_0(z))\Lambda^T
K_{m+1}(c_l,c_l)]^{-1}\pmb{\delta}_z ({\mathcal{P}}_0(z))\Lambda^T.$$
Then
$$\widetilde{V}_m(h(x))= D_m V_m(h(x))-L_mK_{m+1}(x,c_l).$$
Using the definition of the kernel as well as the three-term recurrence relation that  $\{V_m\}_{m \in {\mathbb{N}}}$ satisfies, then
\begin{multline*} z\widetilde{V}_m(z)=
[D_m A_m -L_m G_m(0) A_m ]V_{m+1}(z)  \\+ [D_m B_m +L_m
G_{m+1}(0)C_{m+1}]V_m(z)+D_m C_m V_{m-1}(z).
\end{multline*}
Then the  comparison of the coefficients in~\eqref{lixo} leads to the representation of $\alpha^{1}_{m+1}$, $\alpha^{2}_{m},$ and $\alpha^{3}_{m-1}$.
 \end{proof}
 \begin{example}
In the last theorem if we consider the generalized second kind Chebyshev polynomials, $\{U_m^{A,B,C}\}_{m \in {\mathbb{N}}}$, left-ortho\-gonal with respect to the generalized Markov function ${\mathcal{F}}_{A,B,C}$ and satisfying the three-term recurrence relation
$$zU_m^{A,B,C}(z)=A U_{m+1}^{A,B,C}(z) + B U_m^{A,B,C}(z)+ C U_{m-1}^{A,B,C}(z), \quad m\geq 1,$$
and consider the matrix polynomials $\{\widetilde{U}_m^{A,B,C}\}_{m \in {\mathbb{N}}}$ left-ortho\-gonal with respect to
${\widetilde{\mathcal{F}}}_{A,B,C}$ given by
$$z\widetilde{{\mathcal{F}}}_{A,B,C}(z)=z{\mathcal{F}}_{A,B,C}(z)+
\pmb{\delta}_x ({\mathcal{P}}_0(x))\Lambda^T$$ with
${\mathcal{P}}_0(x) = [1\, x\, \cdots \, x^{N-1} ]^T $ and
$$\pmb{\delta} = [\delta_{c_1}\,
\delta^{\prime}_{c_1}\, \cdots \, \delta^{(M_1)}_{c_1}\, \delta_{c_2}
\, \delta^{\prime}_{c_2}\, \cdots \,\delta^{(M_2)}_{c_2} \, \cdots \,
\delta_{c_M}\, \delta^{\prime}_{c_M} \,\cdots \, \delta^{(M_M)}_{c_M}]^T,$$
we get that
 \begin{eqnarray*}
z\widetilde{U}_m^{A,B,C}(z)= \alpha^{1}_{m+1}U_{m+1}^{A,B,C}(z)
+\alpha^{2}_{m} U_{m}^{A,B,C}(z)+\alpha^{3}_{m-1} U_{m-1}^{A,B,C}(z),
 \end{eqnarray*}
with $\alpha^{1}_{m+1}$, $\alpha^{2}_{m}$,
and $\alpha^{3}_{m-1}$ given~by
\begin{eqnarray*} \alpha^{1}_{m+1}&=&
(\widetilde{\beta}_m^m)^{-1} \beta_m^m A -L_m \beta_0^m A,
\\ \alpha^{2}_{m}&=& (\widetilde{\beta}_m^m)^{-1} \beta_m^m B +L_m
\beta_0^{m+1} C, \\ \alpha^{3}_{m-1}&=& (\widetilde{\beta}_m^m)^{-1} \beta_m^m C ,
\end{eqnarray*}
and
 $$
L_m=(\widetilde{\beta}_m^m)^{-1} \beta_m^m U_{m}^{A,B,C}(0)[I_{N
\times N}+ \pmb{\delta}_z ({\mathcal{P}}_0(z))\Lambda^T
K_{m+1}(c_l,c_l)]^{-1}\pmb{\delta}_z ({\mathcal{P}}_0(z))\Lambda^T
 . $$
\end{example}

\section{Relative Asymptotics} \label{sec:4}

The works about perturbations in the matrix Nevai class (see~\cite{MarcellanYakhlefPinar1,MarcellanYakhlefPinar2}) motivated us
 to study the outer relative asymptotics in the case of the generalized matrix Nevai class.

First, we present how can we reinterpret the model studied by the authors from  a vectorial point of view. Second, the outer relative asymptotics is deduced.

Let $\{\widetilde{G}_n\}_{n \in {\mathbb{N}}}$ and $\{\widetilde{\mathcal{B}}_m\}_{m \in {\mathbb{N}}}$ be bi-orthogonal with respect to $\widetilde{{\mathcal{U}}}$, i.e.,
 \begin{eqnarray}\label{mae}
\widetilde{G}_n^T(h(x))\widetilde{{\mathcal{U}}}(\widetilde{{\mathcal{B}}}_m)=I_{N\times
N}\delta_{n,m},
\end{eqnarray}
where
$\widetilde{{\mathcal{U}}}={\mathcal{U}}+ \Lambda \, \pmb{\delta}$ with $\Lambda$ a numerical matrix of dimension $N \times N$,
 \begin{eqnarray*}
 \pmb{\delta} = [\delta_{c_1}\,
\delta^{\prime}_{c_1}\, \cdots \, \delta^{(M_1)}_{c_1}\, \delta_{c_2}
\, \delta^{\prime}_{c_2}\, \cdots \,\delta^{(M_2)}_{c_2} \, \cdots \,
\delta_{c_M}\, \delta^{\prime}_{c_M} \,\cdots \, \delta^{(M_M)}_{c_M}]^T,
 \end{eqnarray*}
and $h$ is defined by~\eqref{h}.

By a straightforward calculation in~\eqref{mae} we have
 \begin{eqnarray*}
(\widetilde{G}_n^T(h(x))({\mathcal{U}}+\Lambda \,\pmb{\delta}_x ))(\widetilde{{\mathcal{B}}}_m)=(\widetilde{G}_n^T(h(x)){\mathcal{U}})(\widetilde{{\mathcal{B}}}_m)+
((\widetilde{G}_n(h(c_l)))^T \Lambda \,\pmb{\delta}_x )(\widetilde{{\mathcal{B}}}_m),
\end{eqnarray*}
where $c_l$ is any zero of the polynomial $h$.

Since $h(c_l)=0$ and
\begin{eqnarray*}
((\widetilde{G}_n(0))^T \Lambda \,\pmb{\delta}_x )(\widetilde{{\mathcal{B}}}_m)=\widetilde{V}_m(0)\pmb{\delta}_x ({\mathcal{P}}_0)\Lambda^T\widetilde{G}_n(0).
\end{eqnarray*}
Thus
\begin{eqnarray*} (\widetilde{G}_n^T(h(x))({\mathcal{U}}+\Lambda \,\pmb{\delta}_x ))(\widetilde{{\mathcal{B}}}_m)=
(\widetilde{G}_n^T(h(x)){\mathcal{U}})(\widetilde{{\mathcal{B}}}_m)+
\widetilde{V}_m(0)\pmb{\delta}_x ({\mathcal{P}}_0)\Lambda^T\widetilde{G}_n(0).
\end{eqnarray*}
By the bi-orthogonality, we have the following matrix interpretation
 \begin{multline*}
 \frac{1}{2\pi i}\int_\gamma \widetilde{V}_m(z)
\widetilde{{\mathcal{F}}}(z) \widetilde{G}_m(z)dz \\
 = \frac{1}{2\pi
i}\int_\gamma \widetilde{V}_m(z) {\mathcal{F}}(z)
\widetilde{G}_m(z)dz+
\widetilde{V}_m(0)\pmb{\delta}_z ({\mathcal{P}}_0(z))\Lambda^T\widetilde{G}_n(0)
 \end{multline*}
where we get the same model studied by F. Marcell\'{a}n, M. Pi\~{n}ar, and  H. O. Yakhlef when we choose all zeros of the polynomial $h$ equal to $c$, i.e., $h(x)=(x-c)^N$.

Next, we will deduce the outer relative asymptotics  of $\{ \widetilde{V}_m(V_m)^{-1} \}_{m \in {\mathbb{N}}} $
when $\{V_m\}_{m \in {\mathbb{N}}}$ belongs to the generalized matrix Nevai class.
 To prove it,  the following lemma is needed as an auxiliary result.
 \begin{lemma}\label{cccc}
Let ${\mathcal{U}}$ and $\widetilde{{\mathcal{U}}} = {\mathcal{U}}+\Lambda\, \pmb{\delta}$ be quasi-definite vector linear functionals and let ${\mathcal{F}}$ and $\widetilde{{\mathcal{F}}}$ be the generalized Markov functions associated to ${\mathcal{U}}$ and $\widetilde{{\mathcal{U}}}$, respectively. The sequences of matrix polynomials bi-orthogonal to ${\mathcal{F}}$ and $\widetilde{{\mathcal{F}}}$ are denoted by $\{V_m\}_{m \in {\mathbb{N}}}$, $\{G_m\}_{m \in {\mathbb{N}}}$, $\{\widetilde{V}_m\}_{m \in {\mathbb{N}},}$ and $\{\widetilde{G}_m\}_{m \in {\mathbb{N}}}$, respectively. Let also, $\alpha^m_m$, $\widetilde{\alpha}^m_m$ and $\beta^m_m$, $\widetilde{\beta}^m_m$ be the leading coefficients  of $h$ of $V_m$, $\widetilde{V}_m$ and $G_m$, $\widetilde{G}_m$, respectively.
 Then,
\begin{multline}\label{bb}
\left((\beta_m^m)^{-1} \widetilde{\beta}_m^m\right) \left(\widetilde{\alpha}^m_m (\alpha_m^m)^{-1}\right) = I_{N \times N} \\- V_m(0) (I_{N \times N} + \pmb{\delta}_z ({\mathcal{P}}_0(z))\Lambda^T
K_{m+1}(c_l,c_l))^{-1}\pmb{\delta}_z ({\mathcal{P}}_0(z))\Lambda^TG_m(0).\end{multline}
\end{lemma}
\begin{proof}
Taking into account~\eqref{fr}, we have $$\widetilde{V}_m(h(x))=D_m V_m(h(x)) - \widetilde{V}_m(0)\pmb{\delta}_z ({\mathcal{P}}_0(z))\Lambda^T K_{m+1}(x,c_l),$$ where $D_m=(\widetilde{\beta}_m^m)^{-1} \beta_m^m$ (see theorem 8 of~\cite{anamendes}), and thus
 \begin{align*}
&\widetilde{\alpha}^m_m (\alpha_m^m)^{-1}=\frac{1}{2 \pi i} \int \widetilde{V}_m(h(x)) {\mathcal{F}}(h(x)) G_m(h(x))d(h(x)) \\
&\phantom{ola}= \frac{1}{2 \pi i} \int (D_m (V_m(h(x))-V_m(0) \left(I_{N \times N} + \pmb{\delta}_z ({\mathcal{P}}_0(z))\Lambda^T
K_{m+1}(c_l,c_l)\right)^{-1}\\
&\phantom{olololololol}\times \pmb{\delta}_z ({\mathcal{P}}_0(z))\Lambda^T K_{m+1}(x,c_l))){\mathcal{F}}(h(x)) G_m(h(x))d(h(x))
 \end{align*}
This means that
\begin{multline*}
D_m^{-1} \widetilde{\alpha}^m_m (\alpha_m^m)^{-1} \\
=I_{N \times N}-  V_m(0)\left(I_{N \times N} + \pmb{\delta}_z ({\mathcal{P}}_0(z))\Lambda^T
K_{m+1}(c_l,c_l)\right)^{-1}\pmb{\delta}_z ({\mathcal{P}}_0(z))\Lambda^T G_m(0) \end{multline*}
and the statement  follows.
\end{proof}
\begin{teo}\label{resultado}
Let ${\mathcal{U}}$ and $\widetilde{{\mathcal{U}}} = {\mathcal{U}}+\Lambda\, \pmb{\delta}$ be quasi-definite vector linear functionals and let ${\mathcal{F}}$ and $\widetilde{{\mathcal{F}}}$ be the generalized Markov functions associated with  ${\mathcal{U}}$ and $\widetilde{{\mathcal{U}}}$, respectively. The sequences of matrix polynomials bi-ortho\-gonal to ${\mathcal{F}}$ and $\widetilde{{\mathcal{F}}}$ are denoted by $\{V_m\}_{m \in {\mathbb{N}}}$, $\{G_m\}_{m \in {\mathbb{N}}}$, $\{\widetilde{V}_m\}_{m \in {\mathbb{N}},}$ and $\{\widetilde{G}_m\}_{m \in {\mathbb{N}}}$, respectively. Let also  $\alpha^m_m$, $\widetilde{\alpha}^m_m$ and $\beta^m_m$, $\widetilde{\beta}^m_m$ be the leading coefficients  of $V_m$, $\widetilde{V}_m$ and $G_m$, $\widetilde{G}_m$, respectively. \\
Assume that $$\lim_{m \to \infty}A_m =A, \quad \lim_{m \to \infty}B_m =B, \quad \mbox{and} \quad \lim_{m \to \infty}C_m =C,$$
where $\{A_m\}_{m \in {\mathbb{N}}}$, $\{B_m\}_{m \in {\mathbb{N}},}$ and $\{C_m\}_{m \in {\mathbb{N}}}$ are the sequences of numerical matrices involved in the recurrence relation~\eqref{ana3}.
Then,
$$\lim_{m \to \infty} (\beta_m^m)^{-1} \widetilde{\beta}_m^m \widetilde{\alpha}^m_m (\alpha_m^m)^{-1}= I_{N \times N} + {\mathcal{F}}_{A,B,C}(0)({\mathcal{F}}^{\prime}_{A,B,C}(0)^{-1}{\mathcal{F}}_{A,B,C}(0).$$
\end{teo}

To prove this theorem we need the following lemma.
\begin{lemma}\label{tania}
Let $\{V_m\}_{m \in {\mathbb{N}}}$ and $\{G_m\}_{m \in {\mathbb{N}}}$ be a sequences of matrix bi-ortho\-gonal polynomials with respect to the generalized Markov function ${\mathcal{F}}$. Let $A_m$, $B_m,$ and $C_m$ be the matrix coefficients that appear in the recurrence relation~\eqref{ana3} such that
$$\lim_{m \to \infty}A_m =A, \quad \lim_{m \to \infty}B_m =B, \quad \mbox{and} \quad \lim_{m \to \infty}C_m =C,$$
where $A$ and $C$ are non-singular matrices, then
$$\lim_{m \to  \infty} G_{m}^{-1}(z)=0_{N \times N}\quad \mbox{and} \quad \lim_{m \to  \infty} V_{m}^{-1}(z)=0_{N \times N}$$
locally  uniformly  in compact subsets of ${\mathbb{C}}\setminus \Gamma$, where $\Gamma=\cap_{N\geq 0} M_N,\,
M_N=\overline{\cup_{m\geq N}Z_m},$ and  $Z_m$ is the set of the zeros of $V_m$.
\end{lemma}
\begin{proof}
 From Liouville-Ostrogradski type formula we have that
\begin{eqnarray*}
V_m^{-1}A_{m-1}^{-1}G_{m-1}^{-1}&=& V_m^{-1}\left({\mathcal{B}}^{(1)}_{m-1}G_{m-1} - V_{m} G^{(1)}_{m-2}\right)G_{m-1}^{-1} \\
&=& V_m^{-1}{\mathcal{B}}^{(1)}_{m-1}- G^{(1)}_{m-2}G_{m-1}^{-1}.
\end{eqnarray*}
 From the generalized Markov theorem for $V_m$ and its analogue for $G_m$, i.e.,
$$\lim_{m \to  \infty} V^{-1}_m(z) {\mathcal{B}}^{(1)}_{m-1} (z) = {\mathcal{F}}(z), \quad \mbox{and} \quad \lim_{m \to  \infty} G^{(1)}_{m-1}(z) G^{-1}_{m} (z) = {\mathcal{F}}(z)$$
we have that
$$\lim_{m \to  \infty} V_m^{-1}(z)A_{m-1}^{-1}G_{m-1}^{-1}(z)=0_{N \times N}.$$
But
\begin{eqnarray*}
V_m^{-1}(z)A_{m-1}^{-1}G_{m-1}^{-1}(z)=V_{m-1}^{-1}(z)V_{m-1}(z)V_m^{-1}(z)A_{m-1}^{-1}G_{m-1}^{-1}(z)
\end{eqnarray*}
According to  Theorem~\ref{maintheorem}  $$\lim_{m  \to \infty} V_{m-1}(z) V_m^{-1}(z) A_{m-1}^{-1}= {\mathcal{F}}_{A,B,C}(z)$$
and since $V^{-1}_{m}$ is bounded for the spectral norm~$\| . \|_2$ (see Lemma~\ref{kill}), we have that  $\lim_{m \to  \infty} G^{-1}_m(0)=0_{N \times N}$.
To prove  $\lim_{m \to  \infty} V^{-1}_m(0)=0_{N \times N}$ we follow an analogue way.
\end{proof}
 {\noindent}{\sc Proof of Theorem~\ref{resultado}}.
Writing $\Phi(c_l)=(\beta_m^m)^{-1} \widetilde{\beta}_m^m \widetilde{\alpha}^m_m (\alpha_m^m)^{-1}$, by Lemma~\ref{cccc} we have
\begin{multline*}
\Phi_m(c_l) \\
= I_{N \times N}-[ G^{-1}_m(0) ((\Lambda\pmb{\delta}_z )({\mathcal{P}}_0(z)))^{-1}V_m^{-1}(0) +
G^{-1}_m(0)K_{m+1}(c_l,c_l))V_m^{-1}(0)]^{-1}
\end{multline*}
To analyze $\displaystyle \lim_{m \to  \infty} \Phi_m(c_l)$ we start by proving that
\begin{eqnarray}\label{seq1}
\lim_{m \to  \infty}G^{-1}_m(0)K_{m+1}(c_l,c_l))V_m^{-1}(0)=- {\mathcal{F}}^{-1}_{A,B,C}(0){\mathcal{F}}^{\prime}_{A,B,C}(0) {\mathcal{F}}^{-1}_{A,B,C}(0).\end{eqnarray}
If we put $G^{-1}_m(0)K_{m+1}(c_l,c_l))V_m^{-1}(0)=\gamma_m(c_l)$, taking into account~\eqref{cdm11}
and that $V_m^{\prime}(z)V_m^{-1}(z)=-V_m(z)(V_m^{-1}(z))^{\prime}$, then we have
\begin{eqnarray*}
\gamma_m(c_l)&=& A_m V^{\prime}_{m+1}(0)V_m^{-1}(0) + G_m^{-1}(0)G_{m+1}(0)C_{m+1}V_m(0)(V_m^{-1}(0))^{\prime}\\
&=&A_m  V^{\prime}_{m+1}(0)V_m^{-1}(0) + G_m^{-1}(0)G_m(0)A_mV_{m+1}(0)(V_m^{-1}(0))^{\prime}\\
&=& -A_m [(V_m(0)V_{m+1}^{-1}(0))^{-1}]^{\prime}\\
&=&-A_m [(V_m(0)V_{m+1}^{-1}(0))^{-1}(V_m(0)V_{m+1}^{-1}(0))^{-1})^{\prime}(V_m(0)V_{m+1}^{-1}(0))^{-1}]
\end{eqnarray*}
Using~\eqref{assrelativa} and~\eqref{convderivadas}, we have~\eqref{seq1}.

Finally, we just have proved  that $$ \lim_{m \to  \infty} G^{-1}_m(0) ((\Lambda\pmb{\delta}_z )({\mathcal{P}}_0(z)))^{-1}V_m^{-1}(0)=0_{N \times N}.$$
 From Lemma~\ref{tania} we have that $\lim_{m \to  \infty} G^{-1}_m(0)=0_{N \times N}$,  $\lim_{m \to  \infty} V^{-1}_m(0)=0_{N \times N}$,
 and
 \begin{multline*}
\|G^{-1}_m(0) ((\Lambda\pmb{\delta}_z )({\mathcal{P}}_0(z)))^{-1}V_m^{-1}(0)\|_2  \\
\leq \|G^{-1}_m(0)\|_2\|((\Lambda\pmb{\delta}_z )({\mathcal{P}}_0(z)))^{-1}\|_2 \|V_m^{-1}(0)\|_2 \, .
 \end{multline*}
Thus, the result follows.
 \begin{rem}
Notice that in the generalized matrix Nevai class
$$\lim_{m  \to  \infty} (\beta_m^m)^{-1} \widetilde{\beta}_m^m \widetilde{\alpha}^m_m (\alpha_m^m)^{-1}$$
is well  determined and
$$\lim_{m  \to  \infty} (\beta_m^m)^{-1} \widetilde{\beta}_m^m \widetilde{\alpha}^m_m (\alpha_m^m)^{-1}= I_{N \times N} + {\mathcal{F}}_{A,B,C}(0)({\mathcal{F}}^{\prime}_{A,B,C}(0))^{-1}{\mathcal{F}}_{A,B,C}(0).$$
If
 \begin{eqnarray}\label{xi}
\Xi= I_{N \times N} + {\mathcal{F}}_{A,B,C}(0)({\mathcal{F}}^{\prime}_{A,B,C}(0))^{-1}{\mathcal{F}}_{A,B,C}(0)
 \end{eqnarray}
and taking into account  that
$$\lim_{m \to  \infty} \beta_m^m = \lim_{m \to  \infty} (\Delta_m)^{-1} = \lim_{m \to  \infty}(C_m C_{m-1} \cdots C_1 \Delta_0)^{-1}= (C \cdots C \Delta_0)^{-1}, $$
$$\lim_{m \to  \infty} \alpha_m^m = \lim_{m \to  \infty} (\Theta_m)^{-1} = \lim_{m \to  \infty}( \Theta_0 A_{-1} A_{0} \cdots A_m)^{-1}= (\Theta_0 A \cdots A)^{-1},$$
we have that
\begin{eqnarray*}
 \lim_{m \to  \infty} \widetilde{\beta}_m^m \widetilde{\alpha}_m^m   &=& \lim_{m \to  \infty} \beta_m^m \, \Xi \, \lim_{m  \to  \infty} \alpha_m^m = (C \cdots C\Delta_0)^{-1} \Xi (\Theta_0 A \cdots A)^{-1}.
\end{eqnarray*}
 \end{rem}

Now, we are ready to prove the following asymptotic result. So, we set that $ \displaystyle \lim_{m  \to  \infty} (\beta_m^m)^{-1} \widetilde{\beta}_m^m = \Psi$.
\begin{teo}
Let $\{V_m\}_{m \in {\mathbb{N}}}$ be  a sequence of matrix orthogonal polynomials with res\-pect to the generalized Markov function ${\mathcal{F}}$ and $\{A_m\}_{m \in {\mathbb{N}}}$, $\{B_m\}_{m \in {\mathbb{N}},}$ and $\{C_m\}_{m \in {\mathbb{N}}}$ the sequences of numerical matrices involved  in the recurrence relation~\eqref{ana3} such that
$$\lim_{m  \to  \infty}A_m =A, \quad \lim_{m  \to  \infty}B_m =B, \quad \mbox{and} \quad \lim_{m  \to  \infty}C_m =C.$$
Then, there exists a sequence of matrix polynomials $\{\widetilde{V}_m\}_{m \in {\mathbb{N}}}$ left-ortho\-gonal with respect to $\widetilde{{\mathcal{F}}}$ defined by~\eqref{ftilde} such that, for ${\mathbb{C}} \setminus \{\Gamma \cup \{c_{1}, c_{2}, \ldots , c_{M}\}\},$ we~have
\begin{multline}\label{fim}
\lim_{m  \to  \infty} \widetilde{V}_m(z) V_m^{-1}(z) \\
 = \Psi^{-1}[I_{N \times N}-\frac{1}{z}I_{N \times N}- \Xi ({\mathcal{F}}^{-1}_{A,B,C} (0) -  ({\mathcal{F}}_{A,B,C} (z))^{-1})],\end{multline}
where $\Xi$ is given by~\eqref{xi}.\end{teo}
\begin{proof}
If $$I_{N \times N} + \pmb{\delta}_z ({\mathcal{P}}_0(z))\Lambda^T
K_{m+1}(c_l,c_l), \,\, m \in {\mathbb{N}},$$ is non-singular for every $c_l,$ a zero of $h$,  then
$$
\widetilde{V}_m(h(x))
 =
D_m V_m(h(x)) - \widetilde{V}_m(0) \pmb{\delta}_z ({\mathcal{P}}_0(z)) \Lambda^T K_{m+1}(x,c_l) , $$
where
 $$\widetilde{V}_m(0) = D_m V_m(0) (I_{N \times N} + \pmb{\delta}_z ({\mathcal{P}}_0(z))\Lambda^T
K_{m+1}(c_l,c_l) )^{-1} $$
 and~$D_m = (\widetilde{\beta}_m^m)^{-1} \beta_m^m$.
Then, multiplying the last relation on the right by $V_m^{-1}$ we obtain
 \begin{multline*}
\widetilde{V}_m(h(x))V_m^{-1}(h(x))=D_m (I_{N \times N}-V_m(0) \left(I_{N \times N} \right. \\  \left. +  \pmb{\delta}_z ({\mathcal{P}}_0(z))\Lambda^T
K_{m+1}(c_l,c_l)\right)^{-1} \pmb{\delta}_z ({\mathcal{P}}_0(z))\Lambda^T K_{m+1}(x,c_l))V_m^{-1}(h(x))
\end{multline*}
But, from the Christoffel-Darboux formula~\eqref{cdm1}
$$K_{m+1}(x,c_l)=\left[G_m(0)A_m
V_{m+1}(h(x))-G_{m+1}(0) C_{m+1}V_m(h(x))\right]/{h(x)},$$
and then, $\widetilde{V}_m V_m^{-1}$ can be written as
\begin{multline*}
\widetilde{V}_m(h(x))V_m^{-1}(h(x))=D_m (I_{N \times N}-\frac{V_m(0)}{h(x)} \\ \times \left(I_{N \times N} + \pmb{\delta}_z ({\mathcal{P}}_0(z))\Lambda^T
K_{m+1}(c_l,c_l)\right)^{-1}\pmb{\delta}_z ({\mathcal{P}}_0(z))\Lambda^TG_m(0)\\ \times (A_m V_{m+1}(h(x))V_m^{-1}(h(x)) - G_m^{-1}(0) G_{m+1}(0) C_{m+1}).
\end{multline*}
Using~\eqref{bb} in the last relation, we get
\begin{multline}\label{bbb}
\widetilde{V}_m(z)V_m^{-1}(z)= [(\beta_m^m)^{-1} \widetilde{\beta}_m^m]^{-1} (I_{N \times N}-\frac{1}{z}I_{N \times N}\\+[(\beta_m^m)^{-1} \widetilde{\beta}_m^m] [\widetilde{\alpha}^m_m (\alpha_m^m)^{-1}] (A_m V_{m+1}(z)V_m^{-1}(z) - G_m^{-1}(0) G_{m+1}(0) C_{m+1}).
\end{multline}
If we denote
\begin{eqnarray*}
\pmb{E}_m(z) &=& G_m^{-1}(0) G_{m+1}(0) C_{m+1} - A_m V_{m+1}(z))V_m^{-1}(z)\\
&=& C^{-1}_{m+1} G_{m+1}^{-1}(0) G_{m}(0) - V_m(z) V_{m+1}^{-1}(z) A_m^{-1}
\end{eqnarray*}
then~\eqref{bbb} becomes
\begin{multline*}
\widetilde{V}_m(z)V_m^{-1}(z) \\
 = [ (\beta_m^m)^{-1} \widetilde{\beta}_m^m]^{-1} ( I_{N \times N}-\frac{1}{z}I_{N \times N} - [ (\beta_m^m)^{-1} \widetilde{\beta}_m^m ] [ \widetilde{\alpha}^m_m (\alpha_m^m)^{-1} ] \pmb{E}_m(z) ).
\end{multline*}
 From Theorem~\ref{maintheorem} and its analog for the sequence of matrix polynomials $\{G_m\}_{m \in {\mathbb{N}}},$  we have
$$\lim_{m  \to  \infty} \pmb{E}_m(z) =({\mathcal{F}}_{A,B,C} (0))^{-1} -  ({\mathcal{F}}_{A,B,C} (z))^{-1}.$$
Since, there exists a sequence of matrix orthogonal polynomials $\{\widetilde{V}_m\}_{m \in {\mathbb{N}}}$ such that
$\displaystyle \lim_{m  \to  \infty} (\beta_m^m)^{-1} \widetilde{\beta}_m^m = \Psi,$ we obtain~\eqref{fim}.
\end{proof}

\ifx\undefined\bysame
\newcommand{\bysame}{\leavevmode\hbox to3em{\hrulefill}\,}
\fi

\end{document}